\documentclass[12pt]{article}
\usepackage{amsmath}
\usepackage{amssymb}
\usepackage{latexsym,bm}
\usepackage{amsthm}
\usepackage{graphicx}

\usepackage[top=2.4cm,bottom=2cm,left=2.4cm,right=2.4cm]{geometry}
\usepackage{algorithm}
\usepackage{algorithmic}
\usepackage{cite}
\allowdisplaybreaks
\usepackage{multirow}
\usepackage{threeparttable}

\usepackage[figuresright]{rotating}
\usepackage{mathrsfs}

\UseRawInputEncoding

\providecommand{\U}[1]{\protect\rule{.1in}{.1in}}
\providecommand{\U}[1]{\protect\rule{.1in}{.1in}}
\providecommand{\U}[1]{\protect\rule{.1in}{.1in}}
\providecommand{\U}[1]{\protect\rule{.1in}{.1in}}

\newtheorem{theorem}{Theorem}[section]

\theoremstyle{definition}
\newtheorem{definition}{Definition}[section]
\theoremstyle{remark}
\newtheorem{remark}{Remark}[section]
\numberwithin{equation}{section}
\theoremstyle{example}

\numberwithin{equation}{section}

\begin{document}

\title{Convergence analysis of a relaxed inertial alternating minimization algorithm with applications}

\author{ Yang Yang$^{a}$, Yuchao Tang$^{a}$, Jigen Peng$^{b}$ \\
\small ${^a}$ Department of Mathematics, Nanchang University, \\
\small Nanchang 330031, P.R. China \\
\small ${^b}$ School of Mathematics and Information Science, Guangzhou University, \\
\small Guangzhou 510006, P.R. China
}

 \date{}

\maketitle

{}\textbf{Abstract.}
The alternating direction method of multipliers (ADMM) is a popular method for solving convex separable minimization problems with linear equality constraints. The generalization of the two-block ADMM to the three-block ADMM is not trivial since the three-block ADMM is not convergence in general. Many variants of three-block ADMM have been developed with guarantee convergence. Besides the ADMM, the alternating minimization algorithm (AMA) is also an important algorithm for solving the convex separable minimization problem with linear equality constraints. The AMA is first proposed by Tseng, and it is equivalent to the forward-backward splitting algorithm applied to the corresponding dual problem. In this paper, we design a variant of three-block AMA, which is derived by employing an inertial extension of the three-operator splitting algorithm to the dual problem. Compared with three-block ADMM, the first subproblem of the proposed algorithm only minimizing the Lagrangian function. As a by-product, we obtain a relaxed algorithm of Davis and Yin. Under mild conditions on the parameters, we establish the convergence of the proposed algorithm in infinite-dimensional Hilbert spaces. Finally, we conduct numerical experiments on the stable principal component pursuit (SPCP) to verify the efficiency and effectiveness of the proposed algorithm.



\section{Introduction}

Many problems in signal and image processing can be modeled as convex minimization problems, whose objective functions may be two-block separable with linear equality constraints. The alternating direction method of multipliers (ADMM) dated back to the work of Glowinski et al. \cite{Glowinski1975} and Gabay et al. \cite{Gabay1976} is a widely used method for solving two-block separable convex minimization problems with linear equality constraints. The ADMM was received much attention in recent years due to its simplicity in solving various inverse problems arising in image restoration and medical image reconstruction. See for example \cite{wangyl2008,YangSIAM2009,wucl_SIAM_2010}. We refer interested readers to \cite{He2012SIAMNA,Monteiro2013SIAMJO,Fang2015MPC,He2015NM} for theoretical results on ADMM with two-block including convergence analysis and convergence rates analysis.

Since the popularity of the two-block ADMM, it is natural to consider how to generalize it to solve a three-block separable convex minimization problem. There exist many problems that suitable for representing in the formulation of three-block other than two-block. For instance, the stable principal component pursuit \cite{Zhou2010ISIT}, the latent variable Gaussian graphical model selection \cite{Chandrasekaran2012AOS}, the robust principal component analysis model with noisy and incomplete data \cite{Tao2011SIAMJOO}, and so on. The three-block separable convex minimization problem is modeled as follows:
\begin{equation}\label{problem1}
\begin{aligned}
 \min_{x_{1},x_{2},x_{3}}\, & f_{1}(x_{1}) + f_{2}(x_{2}) + f_{3}(x_{3})  \\
\textrm{s.t.}\, & L_{1}x_{1}+L_{2}x_{2}+L_{3}x_{3}=b,
\end{aligned}
\end{equation}
where $f_{i}: H_{i}\rightarrow (-\infty,+\infty]$ with $i=1,2,3$ are proper, lower semi-continuous convex functions (not necessarily smooth); $L_{i}: H_{i}\rightarrow H$ with $i=1,2,3$ are bounded linear operators; and $b\in H$ is a given vector; $H$ and $H_{i}$ with $i=1,2,3$ are real Hilbert spaces. Throughout this paper, we assume that the solution set of problem (\ref{problem1}) exists. For solving the convex minimization problem (\ref{problem1}), the direct extension of the three-block ADMM iterative scheme is as follows:
\begin{equation}\label{3-block-ADMM}
   \left\{
\begin{aligned}
x_{1}^{k+1} &  = \arg\min_{x_{1}} \{ f_{1}(x_{1})  + \frac{\gamma}{2} \| L_{1}x_{1} + L_{2}x_{2}^k + L_{3}x_{3}^k -b -\frac{1}{\gamma}w^{k}\|^2  \}, \\
x_{2}^{k+1} & = \arg\min_{x_{2}} \{ f_{2}(x_{2})  + \frac{\gamma}{2} \| L_{1}x_{1}^{k+1} + L_{2}x_{2} + L_{3}x_{3}^k -b -\frac{1}{\gamma}w^{k} \|^2  \}, \\
x_{3}^{k+1} & = \arg\min_{x_{3}} \{ f_{3}(x_{3})  + \frac{\gamma}{2} \| L_{1}x_{1}^{k+1} + L_{2}x_{2}^{k+1} + L_{3}x_{3} -b -\frac{1}{\gamma}w^{k} \|^2  \}, \\
w^{k+1} & = w^k - \gamma(L_{1}x_{1}^{k+1} + L_{2}x_{2}^{k+1} + L_{3}x_{3}^{k+1} -b),
\end{aligned}
\right.
\end{equation}
where $w$ is the Lagrange multiplier and $\gamma>0$ is the penalty parameter. However, Chen et al. \cite{Chench2016MP} showed that the direct extension of the two-block ADMM to three-block ADMM is divergent if no further condition is imposed. Therefore, many efforts have been made to overcome this shortage. We can roughly divide them into two categories. The first is to make some minor changes to the direct extension of the three-block ADMM. For example, in \cite{He2012SJOO,He2012Manuscript}, He et al. generated a new iteration point by correcting the output of each step to guarantee the convergence of three-block ADMM. He et al. \cite{He2012SJOO} proposed an alternating direction method for prediction correction (ADM-G), which guarantees convergence by adding a Gaussian back substitution correction step. ADM-G specific iteration format is read as:
\begin{equation}\label{3-block-ADM-G}
   \left\{
\begin{aligned}
\tilde{x_{1}}^{k} &  = \arg\min_{x_{1}} \{ f_{1}(x_{1})  + \frac{\gamma}{2} \| L_{1}x_{1} + L_{2}x_{2}^k + L_{3}x_{3}^k -b -\frac{1}{\gamma}w^{k}\|^2  \}, \\
\tilde{x_{2}}^{k} & = \arg\min_{x_{2}} \{ f_{2}(x_{2}) + \frac{\gamma}{2} \| L_{1}\tilde{x_{1}}^{k} + L_{2}x_{2} + L_{3}x_{3}^k -b-\frac{1}{\gamma}w^{k} \|^2  \}, \\
\tilde{x_{3}}^{k} & = \arg\min_{x_{3}} \{ f_{3}(x_{3})  + \frac{\gamma}{2} \| L_{1}\tilde{x_{1}}^{k} + L_{2}\tilde{x_{2}}^{k} + L_{3}x_{3} -b -\frac{1}{\gamma}w^{k}\|^2  \}, \\
\tilde{w}^{k} & = w^k - \gamma(L_{1}\tilde{x_{1}}^{k} + L_{2}\tilde{x_{2}}^{k} + L_{3}\tilde{x_{3}}^{k} -b),\\
x_{1}^{k+1} & = \tilde{x_{1}}^{k},\\
v^{k+1} & = v^{k} - \theta G^{-1}(v^{k}-\tilde{v}^{k}),
\end{aligned}
\right.
\end{equation}
where $\theta \in(0,1)$ and
\begin{equation}\label{def G}
v^{k}=\left(
      \begin{array}{c}
        x_{2}^{k} \\
        x_{3}^{k} \\
        w^{k} \\
      \end{array}
    \right),
    \tilde{v}^{k}=\left(
                \begin{array}{c}
                  \tilde{x_{2}}^{k} \\
                  \tilde{x_{3}}^{k} \\
                  \tilde{w}^{k} \\
                \end{array}
              \right),
    G = \left(
          \begin{array}{ccc}
            I_{2} & (L_{2}^{T}L_{2})^{-1}L_{2}^{T}L_{3} & 0 \\
            0 & I_{3} & 0 \\
            0 & 0 & I \\
          \end{array}
        \right).
\end{equation}
We can see that $G$ is an upper triangular matrix, so step 6 in (\ref{3-block-ADM-G}) is easy to perform. In \cite{Hong2017MP}, Hong and Luo added a contraction factor to the Lagrange multiplier update step and established its global linear convergence under some assumptions. Deng et al. \cite{Deng2017JSC} and Sun et al. \cite{Sun2015SIAMJO} each proposed a variant of three-block ADMM and proved the convergence. Their variants not only add the contraction factor in the Lagrange multiplier update step but also employ an appropriate proximal term in the subproblem of ADMM. The three-block ADMM variant proposed by Sun et al. \cite{Sun2015SIAMJO} uses the Gauss-Seidel cycle to update variables. The three-block ADMM variant proposed by Deng et al. \cite{Deng2017JSC} uses the Jacobi cycle to update variables and considers the general m-block case for any m larger than or equal to 3. The algorithms in \cite{He2012SJOO,He2012Manuscript} belong to the algorithmic framework of prediction-correction methods. We refer interested readers to \cite{Han2015JMIV,Wang2017OMS,Chang2018JCAM,Sun2018JIA,Shen2020Optim} for other types of prediction-correction three-block ADMM. The second way is to add more conditions to the objective function or/and linear equality constraints to ensure the convergence of three-block ADMM. For instance, Han and Yuan in \cite{Han2012JOOTAA} have proved the convergence of three-block ADMM by assuming that the objective functions are strongly convex and the penalty parameter has a small upper bound. In a few years, this condition has been relaxed. The authors of \cite{Chen2013AAA,Lin2015JOTORSOC} proved that convergence of a three-block ADMM iteration scheme if only two of the objective functions are strongly convex and the penalty parameter is limited to a small range. Furthermore, Lin et al. \cite{Lin2014SIAMJOO} proved the globally linear convergence rate of the method under some additional conditions. This paper mainly studies the three-block convex optimization problem with a strongly convex function in the objective function. Cai et al. \cite{Cai2017COA} prove the convergence of (\ref{3-block-ADMM}) when $f_{3}$ is strongly convex with a constant $\mu_{3}>0$, and $L_{1}$, $L_{2}$ are full column rank. Li et al. \cite{LI2015APJOOR} proposed a semi-proximal alternating direction method of multipliers (sPADMM) by hiring appropriate proximal terms on the subproblem of (\ref{3-block-ADMM}) and proved its global convergence. sPADMM iteration details are presented below.
\begin{equation}\label{3-block-sPADMM}
   \left\{
\begin{aligned}
x_{1}^{k+1} &  = \arg\min_{x_{1}} \{ f_{1}(x_{1})  + \frac{\gamma}{2} \| L_{1}x_{1} + L_{2}x_{2}^k + L_{3}x_{3}^k -b -\frac{1}{\gamma}w^{k}\|^2 +\frac{1}{2}\|x_{1}-x_{1}^{k}\|^{2}_{T_{1}} \}, \\
x_{2}^{k+1} & = \arg\min_{x_{2}} \{ f_{2}(x_{2}) + \frac{\gamma}{2} \| L_{1}x_{1}^{k+1} + L_{2}x_{2} + L_{3}x_{3}^k -b -\frac{1}{\gamma}w^{k}\|^2 +\frac{1}{2}\|x_{2}-x_{2}^{k}\|^{2}_{T_{2}} \}, \\
x_{3}^{k+1} & = \arg\min_{x_{3}} \{ f_{3}(x_{3}) + \frac{\gamma}{2} \| L_{1}x_{1}^{k+1} + L_{2}x_{2}^{k+1} + L_{3}x_{3} -b -\frac{1}{\gamma}w^{k}\|^2 +\frac{1}{2}\|x_{3}-x_{3}^{k}\|^{2}_{T_{3}} \}, \\
w^{k+1} & = w^k - \tau\gamma(L_{1}x_{1}^{k+1} + L_{2}x_{2}^{k+1} + L_{3}x_{3}^{k+1} -b),
\end{aligned}
\right.
\end{equation}
where $\tau\in(0,\frac{1+\sqrt{5}}{2})$, $\gamma\in(0,+\infty)$, $T_{i}$ with $i=1,2,3$ are self adjoint and positive semi-definite operators. And function $f_{2}$ is strongly convex with constant $\mu_{2}>0$. The operator $T_{i}$ with $i=1,2,3$ may be $0$ if $\gamma$ is smaller than a threshold. This makes sPADMM (\ref{3-block-sPADMM}) return to  directly extended three-block ADMM (\ref{3-block-ADMM}) with $\tau\in(0,\frac{1+\sqrt{5}}{2})$. The just mentioned three-block ADMM convergence guarantee with strong convexity requirements requires that the penalty parameters be relatively small. Recently, Lin et al. \cite{Lin2017JSC} proved that the convergence of the three-block ADMM has only one strong convexity and smoothness in the objective function and the penalty parameter is larger than zero.

Besides the ADMM and its variants, the alternating minimization algorithm (AMA) proposed by Tseng \cite{Tseng1991SIAM} is an important algorithm for solving a two-block separable convex minimization problem with linear equality constraints, where one of the convex function is assumed to be strongly convex. It is worth noting that the AMA algorithm is equivalent to the forward-backward splitting algorithm applied to the corresponding dual problem. Recently, Davis and Yin \cite{davis2015} proposed a so-called three-block ADMM for solving a three-block separable convex minimization problem, where one of them is strongly convex. The three-block ADMM is derived from the three-operator splitting algorithm applied to the dual problem. They pointed out the three-block ADMM included the Tseng's AMA algorithm, the classical ADMM, and the augmented Lagrangian method. In comparison with the three-block extension of ADMM, the first step of Davis and Yin's three-block ADMM \cite{davis2015} does not involve a quadratic penalty term, which is the same as the AMA algorithm. Therefore, we think it is better to name the three-block ADMM proposed by Davis and Yin \cite{davis2015} as a three-block AMA algorithm. The iteration scheme of the three-block AMA algorithm is as follows:
\begin{equation}\label{3-block-AMA}
   \left\{
\begin{aligned}
x_{1}^{k+1} &  = \arg\min_{x_{1}} \{ f_{1}(x_{1}) - \langle w^k , L_{1}x_{1} \rangle   \}, \\
x_{2}^{k+1} & = \arg\min_{x_{2}} \{ f_{2}(x_{2}) - \langle w^k , L_{2}x_{2} \rangle + \frac{\gamma}{2} \| L_{1}x_{1}^{k+1} + L_{2}x_{2} + L_{3}x_{3}^k -b \|^2  \}, \\
x_{3}^{k+1} & = \arg\min_{x_{3}} \{ f_{3}(x_{3}) - \langle w^k , L_{3}x_{3} \rangle + \frac{\gamma}{2} \| L_{1}x_{1}^{k+1} + L_{2}x_{2}^{k+1} + L_{3}x_{3} -b \|^2  \}, \\
w^{k+1} & = w^k - \gamma(L_{1}x_{1}^{k+1} + L_{2}x_{2}^{k+1} + L_{3}x_{3}^{k+1} -b),
\end{aligned}
\right.
\end{equation}
where $f_{1}$ is a strongly convex function with constant $\mu_{1}>0$. As an important method, the alternating minimization algorithm (AMA) has received extensive attention from scholars.

As Goldstein et al. \cite{goldstein2014SIAM} pointed out that the ADMM and the AMA are preferred ways to solve two-block separable convex programming because of their simplicity, they often perform poorly in situations where the problem is poorly conditioned or when high precision is required. Eckstein and Bertsekas \cite{Eckstein1992} first proposed a relaxed ADMM (RADMM), which included the classical ADMM as a special case. Numerical experiments have been confirmed that the RADMM can accelerate the classical ADMM when the relaxation parameter belongs to $(1,2)$. Further, Xu et al. \cite{xu2017adaptive} proposed an adaptive relaxed ADMM that automatically tuned the algorithm parameters. Goldstein et al. \cite{goldstein2014SIAM} proposed two accelerated variants of the ADMM and the AMA, which are based on Nesterov's accelerated gradient method. Kadkhodaie et al. \cite{Kadkhodaie2015} proposed a so-called accelerated alternating direction method of multipliers (A2DM2) and proved that the algorithm achieved $O(1/k^{2})$ convergence rate, where $k$ is the iteration number. They weakened the assumptions required in \cite{goldstein2014SIAM}. We would like to point out that there also exist some other approaches for accelerating the ADMM, such as accelerated ADMM \cite{Pejcic:227078} based on the Douglas-Rachford envelope (DRE), adaptive accelerated ADMM \cite{Liang2019NIPs}, and accelerated ADMM based on accelerated proximal point algorithm \cite{Kim2019}. Similar to the idea of Nesterov's accelerated gradient method, the inertial method becomes popular in recent years. It provides a general way to select the inertia parameters. Chen et al. \cite{Chen2015} proposed an inertial proximal ADMM, which derived from the inertial proximal point algorithm. On the other hand, Bo\c{t} and Csetnek \cite{Bot2016MTA} proposed an inertial ADMM, which is based on the inertial Douglas-Rachford splitting algorithm \cite{Bot2015AMC}. To the best of our knowledge, we have not seen any generalization work on the three-block AMA algorithm to the relaxation or the inertia.

The purpose of this paper is to introduce a relaxed inertial three-block AMA algorithm for solving the three-block separable convex minimization problem. The idea is to employ the inertial three-operator splitting algorithm \cite{Cui2019} to the dual problem. As a by-product, we obtain a relaxed three-block AMA algorithm, which generalizes the three-block AMA algorithm of Davis and Yin \cite{davis2015}. Under mild conditions, we prove the convergence of the proposed algorithms. To verify the efficiency and effectiveness of the proposed algorithms, we apply them to solve the stable principal component pursuit (SPCP) \cite{Zhou2010ISIT} problem. We also report numerical results compared with other algorithms for solving the SPCP.

We highlight the contributions of this paper: (i) We propose a generalization of the three-block AMA with relaxation and inertia. The obtained algorithm includes several algorithms as its special cases; (ii) We study the convergence of the proposed algorithm under different conditions on the parameters in infinite-dimensional Hilbert spaces. Compared with other existing three-block ADMM and its variants, we obtain weak and strong convergence of the iteration schemes; (iii) We conduct extensive numerical experiments on SPCP to verify the impact of the introduced relaxation and inertia parameters.

The rest of this paper is organized as follows. In Section 2, we present some preliminaries on the maximally monotone operators and convex functions. In particular, we review several results on the inertial three-operator splitting algorithm. In Section 3, we present the main algorithm and prove the convergence of it. In Section 4, we conduct numerical experiments on the stable principal component pursuit to demonstrate the efficiency and effectiveness of the proposed algorithms. Finally, we give some conclusions. We also present two open questions for further study.

\section{Preliminaries}

In this section, we review some basic definitions in convex analysis and monotone operator theory. Most of the following definitions can be found in \cite{bauschkebook2017}. Let $H$ is a real Hilbert space, which endowed with an inner product $\langle\cdot,\cdot\rangle$ and associated norm $\|\cdot\|=\sqrt{\langle\cdot,\cdot\rangle}$. The symbols $\rightharpoonup$ and $\rightarrow$ denote weak and strong convergence, respectively.

Let $A:H\rightarrow2^{H}$ be a set-valued operator. We denote by $zerA=\{x\in H:0\in Ax\}$ its set of zeros, by $graA=\{(x,u)\in H\times H:u\in Ax\}$ its graph and by $ran A=\{u\in H:\exists x\in H,u\in Ax\}$ its range.   The resolvent of an operator $A:H\rightarrow2^{H}$ is denoted by $J_{A}=(I+A)^{-1}$.

\begin{definition}(\cite{bauschkebook2017})
Let $A:H\rightarrow2^{H}$ be a set-valued operator. Then

(i) $A$ is monotone if
$$
\langle x-y,u-v\rangle\geq 0, \forall(x,u),(y,v)\in gra A.
$$
 Further, $A$ is maximally monotone if there exists no monotone operator $A':H\rightarrow 2^{H}$ such thst $graA'$ properly contains $graA$.

(ii) $A$ is uniformly monotone if there exists an increasing function $\phi:[0,+\infty)\rightarrow[0,+\infty]$ that vanishes only at 0 such that
$$
\langle x-y,u-v\rangle\geq\phi(\|x-y\|), \forall(x,u),(y,v)\in graA.
$$

\end{definition}

\begin{definition}(\cite{bauschkebook2017})\label{def2}
Let $B:H\rightarrow H$ is a single valued operator.  $B:H\rightarrow H$ is said to be $\beta-cocoercive$, for some $\beta>0$, if
$$
\langle x-y,Bx-By\rangle \geq \beta\|Bx-By\|^{2}, \forall x,y\in H.
$$

\end{definition}

Let a function $f:H\rightarrow(-\infty,+\infty]$. We denote by $\Gamma_{0}(H)$ be the class of proper, lower semicontinuous convex functions $f:H\rightarrow(-\infty,+\infty]$. Let $f\in\Gamma_{0}(H)$, the conjugate of $f$ is $f^{\ast}\in\Gamma_{0}(H)$ defined by $f^{\ast}:u\mapsto sup_{x\in H}(\langle x,u\rangle-f(x))$, and the subdifferential of $f$ is the maximally monotone operator
$$
\partial f:H\rightarrow2^{H}:x\mapsto \{u\in H|f(y)\geq f(x)+\langle u,y-x\rangle,\forall y\in H\}.
$$
$f$ is uniformly convex if there exists an increasing function $\phi:[0,+\infty)\rightarrow[0,+\infty]$ that vanishes only at 0 such that
$$
\alpha f(x)+(1-\alpha)f(y)\geq f(\alpha x+(1-\alpha)y)+\alpha(1-\alpha)\phi(\|x-y\|), \forall x,y\in H,\alpha\in(0,1).
$$
$f$ is $\sigma-strongly$ convex for some $\sigma>0$ if $f-\frac{\sigma}{2}\|\cdot\|^{2}$ is convex.

The proximity operator $prox_{\lambda f}:x\mapsto \arg\min_{y}\{f(y)+\frac{1}{2\lambda}\|x-y\|^{2}\}$, where $\lambda>0$. Let $f\in\Gamma_{0}(H)$, then we have $J_{\lambda\partial f}=(I+\lambda\partial f)^{-1}=prox_{\lambda f}$.

To analyze the convergence of the algorithm proposed in this paper, we recall the main results of the inertial three-operator splitting algorithm in \cite{Cui2019}.

\begin{theorem}(\cite{Cui2019})\label{theo1}
Let $H$ be real Hilbert space. Let $A, B:H\rightarrow 2^{H}$ be two maximally monotone operators. Let $C:H\rightarrow H$ be a $\beta-cocoercive$ operator, for some $\beta>0$. Let $z^{0}, z^{1}\in H$, and set
\begin{equation}\label{inertial-three-operator}
   \left\{
\begin{aligned}
y^{k}&=z^{k}+\alpha_{k}(z^{k}-z^{k-1}),\\
w^{k}&=J_{\gamma B}y^{k},\\
u^{k}&=J_{\gamma A}(2w^{k}-y^{k}-\gamma Cw^{k}),\\
z^{k+1}&=y^{k}+\lambda_{k}(u^{k}-w^{k}),
\end{aligned}
\right.
\end{equation}
where the parameters $\gamma$, $\{\alpha_{k}\}$ and $\{\lambda_{k}\}$ satisfy the following conditions:

\rm(c1) $\gamma\in(0,2\beta\bar{\varepsilon})$, where $\bar{\varepsilon}\in(0,1)$;

(c2) $\{\alpha_{k}\}$ is nondecreasing with $k\geq1$, $\alpha_{1}=0$ and $0\leq\alpha_{k}\leq\alpha<1$;

(c3) for every $k\geq1$, and $\lambda, \sigma, \delta>0$ such that
\begin{equation}\label{eq2.4}
\delta>\frac{\alpha^{2}(1+\alpha)+\alpha\sigma}{1-\alpha^{2}} \, \textrm{ and } \, 0<\lambda\leq\lambda_{k}\leq\frac{\delta-\alpha[\alpha(1+\alpha)+\alpha\delta+\sigma]}{\bar{\alpha}\delta[1+\alpha(1+\alpha)+\alpha\delta+\sigma]},
\end{equation}
where $\bar{\alpha}=\frac{1}{2-\bar{\varepsilon}}$. Then the following hold:

(i) $\{z^{k}\}$ converges weakly to $z^{\ast}$;

(ii) $\{w^{k}\}$ converges weakly to $J_{\gamma B}z^{\ast}\in zer(A+B+C)$;

(iii) $\{u^{k}\}$ converges weakly to $J_{\gamma B}z^{\ast}=J_{\gamma A}(2J_{\gamma B}z^{\ast}-z^{\ast}-\gamma CJ_{\gamma B}z^{\ast})\in zer(A+B+C)$;

(iv) $\{z^{k}-z^{k-1}\}$ converges strongly to $0$;

(v) $\{w^{k}-u^{k}\}$ converges strongly to $0$;

(vi) $\{Cw^{k}\}$ converges strongly to $CJ_{\gamma B}z^{\ast}$;

(vii) Suppose that one of the following conditions hold:

\quad\quad(a) $A$ be uniformly monotone on every nonempty bounded subset of $dom A$;

\quad\quad(b) $B$ be uniformly monotone on every nonempty bounded subset of $dom B$;

\quad\quad(c) $C$ be demiregular at every point $x\in zer(A+B+C)$.

Then $\{w^{k}\}$ and $\{u^{k}\}$ converge strongly to $J_{\gamma B}z^{\ast}\in zer(A+B+C)$.
\end{theorem}
\vskip 2mm

\begin{proof}
(i), (ii), (iii) and (vii) are directly obtained from Theorem 3.1 of \cite{Cui2019}.
(iv), (v) and (vi) can be easily obtained from Theorem 3.1 of \cite{Cui2019}. Here we omit the proof.
\end{proof}
\vskip 2mm

\begin{theorem}(\cite{Cui2019})\label{theo2}
Let $H$ be a real Hilbert space. Let $A,B:H\rightarrow 2^{H}$ be two maximally monotone operators. Let $C:H\rightarrow H$ be a $\beta-cocoercive$ operator, for some $\beta>0$. Let the iterative sequences $\{z^{k}\}$, $\{w^{k}\}$ and $\{u^{k}\}$ are generated by (\ref{inertial-three-operator}). Assume that the parameters $\gamma$, $\{\alpha_{k}\}$ and $\{\lambda_{k}\}$ satisfy the following conditions:

\rm(c1) $\gamma\in(0,2\beta\bar{\varepsilon})$, where $\bar{\varepsilon}\in(0,1)$;

(c2) $0\leq \alpha_{k}\leq \alpha< 1$ and  $0<\underline{\lambda}\leq \lambda_{k}\bar{\alpha}\leq\overline{\lambda}<1$, where $\bar{\alpha}=\frac{1}{2-\bar{\varepsilon}}$;

(c3) $\sum_{k=0}^{+\infty}\alpha_{k}\|z^{k}-z^{k-1}\|^{2}<+\infty$.

Then the following hold:

\rm(i) $\{z^{k}\}$ converges weakly to $z^{\ast}$;

(ii) $\{w^{k}\}$ converges weakly to $J_{\gamma B}z^{\ast}\in zer(A+B+C)$;

(iii) $\{u^{k}\}$ converges weakly to $J_{\gamma B}z^{\ast}=J_{\gamma A}(2J_{\gamma B}z^{\ast}-z^{\ast}-\gamma CJ_{\gamma B}z^{\ast})\in zer(A+B+C)$;

(iv) $\{\alpha_{k}(z^{k}-z^{k-1})\}$ converges strongly to $0$;

(v) $\{w^{k}-u^{k}\}$ converges strongly to $0$;

(vi) $\{Cw^{k}\}$ converges strongly to $CJ_{\gamma B}z^{\ast}$;

(vii) Suppose that one of the following conditions hold:

\quad\quad(a) $A$ be uniformly monotone on every nonempty bounded subset of $dom A$;

\quad\quad(b) $B$ be uniformly monotone on every nonempty bounded subset of $dom B$;

\quad\quad(c) $C$ be demiregular at every point $x\in zer(A+B+C)$.

Then $\{w^{k}\}$ and $\{u^{k}\}$ converge strongly to $J_{\gamma B}z^{\ast}\in zer(A+B+C)$.
\end{theorem}
\vskip 2mm

\begin{proof}
The proof of Theorem \ref{theo2} is similar to Theorem \ref{theo1}, so we omit it here.
\end{proof}

\begin{remark}\label{rema1}
The condition $\alpha_{1}=0$ in the Theorem \ref{theo1} can be replaced by the assumption $z^{0}=z^{1}$.
\end{remark}
\vskip 2mm

\section{Relaxed inertial three-block AMA for solving three-block separable convex minimization problem }

In this section, we present the main results of this paper including our proposed algorithm and its convergence theorem. The Lagrange function of problem (\ref{problem1}) is defined as follows:
\begin{equation}\label{eq3.1}
L(x_{1},x_{2},x_{3},w)=f_{1}(x_{1}) + f_{2}(x_{2}) + f_{3}(x_{3}) -\langle L_{1}x_{1}+L_{2}x_{2}+L_{3}x_{3}-b,w\rangle,
\end{equation}
where $w$ is a Lagrange multiplier. Through to Lagrange function (\ref{eq3.1}), the dual problem of problem (\ref{problem1}) is
\begin{equation}\label{eq3.2}
\min_{w\in H}f_{1}^{\ast}(L_{1}^{\ast}w)+f_{2}^{\ast}(L_{2}^{\ast}w) +f_{3}^{\ast}(L_{3}^{\ast}w) -\langle b,w\rangle,
\end{equation}
where $f_{i}^{\ast}$   are the Fenchel-conjugate functions of $f_{i}$, respectively. According to the first-order optimality condition of problem (\ref{problem1}), the solution of problem (\ref{problem1}) is equivalent to finding $x_{i}^{\ast}\in H_{i}$ and $w^{\ast}\in H$ satisfying the following formula:
\begin{equation}\label{eq3.3}
\begin{aligned}
 & 0\in \partial f_{1}(x_{1}^{\ast})-L_{1}^{\ast}w^{\ast},  \\
 & 0\in \partial f_{2}(x_{2}^{\ast})-L_{2}^{\ast}w^{\ast},  \\
 & 0\in \partial f_{3}(x_{3}^{\ast})-L_{3}^{\ast}w^{\ast},  \\
 & L_{1}x_{1}^{\ast}+L_{2}x_{2}^{\ast}+L_{3}x_{3}^{\ast}-b=0,
\end{aligned}
\end{equation}
this is what we usually call the KKT condition.

Next, we present the main algorithm of this paper and prove its convergence.

\vskip 2mm
\renewcommand{\algorithmicrequire}{\textbf{Input:}}
\renewcommand{\algorithmicensure}{\textbf{Output:}}
\begin{algorithm}[htb]
\caption{Relaxed inertial three-block AMA }
\label{algo2}
\begin{algorithmic}[1]
\REQUIRE
\vskip 2mm
For arbitrary $w^{1}\in H$, $p^{1}=0$ and $x_{3}^{1}\in H_{3}$. Choose $\gamma$, $\alpha_{k}$ and $\lambda_{k}$.\\
For $k= 1,2,3, \cdots$, compute
\quad  \STATE $x_{1}^{k+1}   = \arg\min_{x_{1}} \{ f_{1}(x_{1})-\langle w^k,L_{1}x_{1}\rangle\}$,

 \STATE $x_{2}^{k+1}  = \arg\min_{x_{2}} \{ f_{2}(x_{2}) - \langle w^k , L_{2}x_{2} \rangle + \frac{\gamma}{2} \| L_{1}x_{1}^{k+1} + L_{2}x_{2} + L_{3}x_{3}^{k} -b \|^2  \}$,

  \STATE $x_{3}^{k+1}  = \arg\min_{x_{3}} \{ f_{3}(x_{3}) -\langle w^k+\alpha_{k+1}p^{k}, L_{3}x_{3} \rangle + \frac{\gamma}{2} \|L_{3}(x_{3}-x_{3}^{k}) +(1+\alpha_{k+1})\lambda_{k}(L_{1}x_{1}^{k+1} + L_{2}x_{2}^{k+1} + L_{3}x_{3}^{k} -b) \|^2\}$,

  \STATE $w^{k+1}  = w^k + \alpha_{k+1}p^{k} - \gamma(L_{3}(x_{3}^{k+1}-x_{3}^{k}) +(1+\alpha_{k+1})\lambda_{k}(L_{1}x_{1}^{k+1} + L_{2}x_{2}^{k+1} + L_{3}x_{3}^{k} -b))$,

  \STATE $p^{k+1} =\alpha_{k+1}(p^{k}-\gamma\lambda_{k}(L_{1}x_{1}^{k+1} + L_{2}x_{2}^{k+1} + L_{3}x_{3}^{k} -b))$.

Stop when a given stopping criterion is met.
\ENSURE $x_{1}^{k+1},x_{2}^{k+1},x_{3}^{k+1}$ and $w^{k+1}$.
\end{algorithmic}
\end{algorithm}
\vskip 2mm

To study the convergence analysis of Algorithm \ref{algo2}, we make the following assumptions:

(A1). Assume that $f_{1}$ is $\mu$-strongly convex, for some $\mu>0$.

(A2). The optimal solution of problem (\ref{problem1}) is nonempty, and the exists $x'=(x_{1}',x_{2}',x_{3}')\in ri(dom(f_{1})\times dom(f_{2})\times dom(f_{3}))\cap C$, where $C=\{(x_{1},x_{2},x_{3})\in H_{1}\times H_{2}\times H_{3}|L_{1}x_{1}+L_{2}x_{2}+L_{3}x_{3}=b\}$.

(A3). For each $i=1,2,3$, let the bounded linear operator $L_{i}$ satisfies that $\|L_{i}x_{i}\|\geq\theta_{i}\|x_{i}\|$, for some $\theta_{i}>0$ and $\forall x_{i}\in H_{i}$.

Under the assumption (A2), we know that the dual solution of problem (\ref{problem1}) is nonempty, and the strong duality holds, i.e., $v(P)=v(D)$.

Next, we will prove the convergence theorem of Algorithm \ref{algo2} under two different conditions. For the convergence proof of Algorithm \ref{algo2}, we roughly divide it into two steps. First, we prove that Algorithm \ref{algo2} is equivalent to (\ref{inertial-three-operator}), and Algorithm \ref{algo2} is derived from (\ref{inertial-three-operator}) through variable substitution. Secondly, after proving the equivalence of Algorithm \ref{algo2} and (\ref{inertial-three-operator}), we can prove the convergence conclusions of the Algorithm \ref{algo2} by using Theorem \ref{theo1} and Theorem \ref{theo2}, respectively.

\begin{theorem}\label{theo3}
Suppose that the assumptions (A1)-(A3) are valid. Let $\{(x_{1}^{k},x_{2}^{k},x_{3}^{k},w^{k})\}$ be the sequence generated by Algorithm \ref{algo2}. Let $\gamma\in(0,2\beta\bar{\varepsilon})$, where $\bar{\varepsilon}\in(0,1)$ and $\beta=\mu/\|L_{1}\|^{2}$. Let $\{\alpha_{k}\}$ is nondecreasing with $\alpha_{1}=0$ and $0\leq\alpha_{k}\leq\alpha<1$. Let $\lambda>0,\sigma>0,\delta>0$ and $\{\lambda_{k}\}$ such that
$$
\delta>\frac{\alpha^{2}(1+\alpha)+\alpha\sigma}{1-\alpha^{2}}\, \textrm{ and } \, 0<\lambda\leq\lambda_{k}\leq\frac{\delta-\alpha[\alpha(1+\alpha)+\alpha\delta+\sigma]}{\bar{\alpha}\delta[1+\alpha(1+\alpha)+\alpha\delta+\sigma]}\quad \forall k\geq1,
$$
where $\bar{\alpha}=\frac{1}{2-\bar{\varepsilon}}$. Then there exists a point pair $(x_{1}^{\ast},x_{2}^{\ast},x_{3}^{\ast},w^{\ast})$, which is the saddle point of the Lagrange function (\ref{eq3.1}) such that the following hold:

\rm(i) $\{(x_{1}^{k+1},x_{2}^{k+1},x_{3}^{k+1})\}_{k\geq1}$ converges weakly to $(x_{1}^{\ast},x_{2}^{\ast},x_{3}^{\ast})$. In particular, $\{x_{1}^{k+1}\}_{k\geq1}$ converges strongly to $x_{1}^{\ast}$;

(ii) $\{w^{k+1}\}_{k\geq1}$ converges weakly to $w^{\ast}$;

(iii) $\{L_{1}x_{1}^{k+1}+L_{2}x_{2}^{k+1}+L_{3}x_{3}^{k}\}_{k\geq2}$ converges strongly to $b$;

(iv) Suppose that one of the following conditions hold:

\quad\quad (a) $f_{1}^{\ast}$ is uniformly convex on every nonempty bounded subset of $dom f_{1}^{\ast}$;

\quad\quad (b) $f_{2}^{\ast}$ is uniformly convex on every nonempty bounded subset of $dom f_{2}^{\ast}$;

\quad\quad (c) $f_{3}^{\ast}$ is uniformly convex on every nonempty bounded subset of $dom f_{3}^{\ast}$;

then $\{w^{k+1}\}_{k\geq1}$ converges strongly to the unique optimal solution of $(D)$;

(v) $\lim_{k\rightarrow+\infty}(f_{1}(x_{1}^{k+1})+f_{2}(x_{2}^{k+1})+f_{3}(x_{3}^{k})) =v(P)=v(D)=\lim_{k\rightarrow+\infty}(-f_{1}^{\ast}(L_{1}^{\ast}w^{k}) -f_{2}^{\ast}(L_{2}^{\ast}u^{k})-f_{3}^{\ast}(L_{3}^{\ast}w^{k}) +\langle w^{k},b\rangle)$, where $u^{k}$ is defined as follows
\begin{equation}\label{eq1}
u^{k}=w^{k}-\gamma(L_{1}x_{1}^{k+1}+L_{2}x_{2}^{k+1}+L_{3}x_{3}^{k}-b).
\end{equation}
\end{theorem}

\vskip 2mm

\begin{proof}
Let $A=\partial(f_{2}^{\ast}\circ(L_{2}^{\ast}\cdot)),\  B=\partial(f_{3}^{\ast}\circ(L_{3}^{\ast}\cdot)-\langle b,\cdot\rangle),\ \textrm{ and } \ C=\nabla(f_{1}^{\ast}\circ(L_{1}^{\ast}\cdot))$. Since $f_{1}$ is $\mu$-strongly convex, $\nabla(f_{1}^{\ast}\circ(L_{1}^{\ast}\cdot))$ is $\mu/\|L_{1}\|^{2}$-cocoercive, and $A,B$ are maximally monotone. Then, we obtain the following inertial three-operator splitting algorithm from \cite{Cui2019} to solve the dual problem (\ref{eq3.2}).
\begin{equation}\label{eq2}
\begin{aligned}
& y^{k}=z^{k}+\alpha_{k}(z^{k}-z^{k-1}), \\
& w^{k}=J_{\gamma B}y^{k}, \\
& u^{k}=J_{\gamma A}(2w^{k}-y^{k}-\gamma Cw^{k}), \\
& z^{k+1}=y^{k}+\lambda_{k}(u^{k}-w^{k}).
\end{aligned}
\end{equation}

Next, we prove that the iterative sequence $\{(x_{1}^{k},x_{2}^{k},x_{3}^{k},w^{k})\}$ generated by Algorithm \ref{algo2} is equivalent to the inertial three-operator splitting algorithm (\ref{eq2}).

From $w^{k}=J_{\gamma B}y^{k}$ and $B=\partial(f_{3}^{\ast}\circ(L_{3}^{\ast}\cdot)-\langle b,\cdot\rangle)$, we have
\begin{equation}\label{eq3}
y^{k}-w^{k}\in\gamma(L_{3}\partial f_{3}^{\ast}(L_{3}^{\ast}w^{k})-b).
\end{equation}
Let $x_{3}^{k}\in\partial f_{3}^{\ast}(L_{3}^{\ast}w^{k})$, we get
\begin{equation}\label{eq4}
L_{3}^{\ast}w^{k}\in\partial f_{3}(x_{3}^{k}),\ \textrm{and} \ z^{k}+\alpha_{k}(z^{k}-z^{k-1})-w^{k}=\gamma L_{3}x_{3}^{k}-\gamma b.
\end{equation}

It follows from the definition of $C$ that $Cw^{k}=L_{1}\nabla f_{1}^{\ast}(L_{1}^{\ast}w^{k})$. Let $x_{1}^{k+1}=\nabla f_{1}^{\ast}(L_{1}^{\ast}w^{k})$, then we obtain
\begin{equation}\label{eq5}
L_{1}^{\ast}w^{k}\in\partial f_{1}(x_{1}^{k+1}),\ \textrm{and} \ Cw^{k}=L_{1}x_{1}^{k+1}.
\end{equation}

From the relation $L_{1}^{\ast}w^{k}\in\partial f_{1}(x_{1}^{k+1})$ yields,
\begin{equation}\label{eq6}
x_{1}^{k+1}=\arg\min_{x_{1}}\{f_{1}(x_{1})-\langle w^{k},L_{1}x_{1}\rangle\},
\end{equation}
which is the first step of Algorithm \ref{algo2}.

Notice that $u^{k}=J_{\gamma A}(2w^{k}-y^{k}-\gamma Cw^{k})$ and (\ref{eq5}), we have
\begin{equation}\label{eq7}
2w^{k}-y^{k}-\gamma L_{1}x_{1}^{k+1}-u^{k}\in\gamma L_{2}\partial f_{2}^{\ast}(L_{2}^{\ast}u^{k}).
\end{equation}
Let $x_{2}^{k+1}\in\partial f_{2}^{\ast}(L_{2}^{\ast}u^{k})$, we obtain that
\begin{equation}\label{eq8}
L_{2}^{\ast}u^{k}\in\partial f_{2}(x_{2}^{k+1})\ \textrm{and} \ 2w^{k}-y^{k}-\gamma L_{1}x_{1}^{k+1}-u^{k}=\gamma L_{2}x_{2}^{k+1}.
\end{equation}

By combining (\ref{eq4}) and (\ref{eq8}), we get
\begin{equation}\label{eq9}
u^{k}=w^{k}-\gamma(L_{1}x_{1}^{k+1}+L_{2}x_{2}^{k+1}+L_{3}x_{3}^{k}-b).
\end{equation}

Consequently, we obtain
\begin{align}\label{eq10}
& 0\in\partial f_{2}(x_{2}^{k+1})-L_{2}^{\ast}(w^{k}- \gamma(L_{1}x_{1}^{k+1}+L_{2}x_{2}^{k+1}+L_{3}x_{3}^{k}-b)),  \nonumber\\
\Leftrightarrow & x_{2}^{k+1}=\arg\min_{x_{2}}\{f_{2}(x_{2})-\langle w^{k},L_{2}x_{2}\rangle+\frac{\gamma}{2}\| L_{1}x_{1}^{k+1}+L_{2}x_{2}+L_{3}x_{3}^{k}-b\|^{2}\},
\end{align}
which is the second step of Algorithm \ref{algo2}.

Set $p^{k}=\alpha_{k}(z^{k}-z^{k-1})$. Then it follows from (\ref{eq4}) that
\begin{equation}\label{eq11}
p^{k}=w^{k}+\gamma L_{3}x_{3}^{k}-\gamma b-z^{k}.
\end{equation}

Further, we have
\begin{align}\label{eq12}
p^{k+1} & =\alpha_{k+1}(z^{k+1}-z^{k})  \nonumber\\
& =\alpha_{k+1}(w^{k+1}+\gamma L_{3}x_{3}^{k+1}-\gamma b-p^{k+1}-(w^{k}+\gamma L_{3}x_{3}^{k}-\gamma b-p^{k})),
\end{align}
which implies that
\begin{equation}\label{eq13}
p^{k+1}=\frac{\alpha_{k+1}}{1+\alpha_{k+1}}(p^{k}+(w^{k+1}-w^{k})+\gamma L_{3}(x_{3}^{k+1}-x_{3}^{k})).
\end{equation}

By (\ref{eq4}) and $z^{k+1}=y^{k}+\lambda_{k}(u^{k}-w^{k})$, we have
\begin{align}\label{eq14}
z^{k}+\alpha_{k}(z^{k}-z^{k-1})+\lambda_{k}(u^{k}-w^{k}) & =z^{k+1}  \nonumber\\
& =w^{k+1}+\gamma L_{3}x_{3}^{k+1}-\gamma b-\alpha_{k+1}(z^{k+1}-z^{k}).
\end{align}

From (\ref{eq4}), (\ref{eq9}), (\ref{eq13}) and (\ref{eq14}), we get
\begin{align}\label{eq15}
w^{k+1} & = z^{k}+\alpha_{k}(z^{k}-z^{k-1})+\lambda_{k}(u^{k}-w^{k}) +\alpha_{k+1}(z^{k+1}-z^{k})-\gamma L_{3}x_{3}^{k+1}+\gamma b  \nonumber\\
& = w^{k}+\gamma L_{3}x_{3}^{k}-\gamma b -\gamma \lambda_{k}(L_{1}x_{1}^{k+1}+L_{2}x_{2}^{k+1}+L_{3}x_{3}^{k}-b) +p^{k+1}-\gamma L_{3}x_{3}^{k+1}+\gamma b.
\end{align}
Consequently, we obtain
\begin{equation}\label{eq16}
w^{k+1}=w^{k}+\alpha_{k+1}p^{k}-\gamma(L_{3}(x_{3}^{k+1}-x_{3}^{k})+ \lambda_{k}(1+\alpha_{k+1})(L_{1}x_{1}^{k+1}+L_{2}x_{2}^{k+1}+L_{3}x_{3}^{k}-b)),
\end{equation}
which is the fourth step of Algorithm \ref{algo2}.

Combining (\ref{eq13}) with (\ref{eq16}), we get
\begin{equation}\label{eq16.5}
p^{k+1} =\alpha_{k+1}(p^{k}-\gamma\lambda_{k}(L_{1}x_{1}^{k+1} + L_{2}x_{2}^{k+1} + L_{3}x_{3}^{k} -b)),
\end{equation}
which is the fifth step of Algorithm \ref{algo2}.

According to $L_{3}^{\ast}w^{k+1}\in\partial f_{3}(x_{3}^{k+1})$, we have
\begin{equation}\label{eq17}
\begin{aligned}
0\in\partial f_{3}&(x_{3}^{k+1}) -L_{3}^{\ast}\{w^{k}+\alpha_{k+1}p^{k} \\
&-\gamma[L_{3}(x_{3}^{k+1}-x_{3}^{k})+ \lambda_{k}(1+\alpha_{k+1})(L_{1}x_{1}^{k+1}+L_{2}x_{2}^{k+1}+L_{3}x_{3}^{k}-b)]\},
\end{aligned}
\end{equation}
which is equivalent to
\begin{equation}\label{eq18}
\begin{aligned}
x_{3}^{k+1}=\arg\min_{x_{3}}&\{ f_{3}(x_{3})-\langle w^{k}+\alpha_{k+1}p^{k},L_{3}x_{3}\rangle \\
+& \frac{\gamma}{2}\|L_{3}(x_{3}-x_{3}^{k})+\lambda_{k}(1+\alpha_{k+1})
(L_{1}x_{1}^{k+1}+L_{2}x_{2}^{k+1}+L_{3}x_{3}^{k}-b)\|^{2}\},
\end{aligned}
\end{equation}
and is the third step of Algorithm \ref{algo2}. Therefore, we can conclude from the above that Algorithm \ref{algo2} is equivalent to (\ref{eq2}).

It follows from Theorem \ref{theo1} that there exists $z^{\ast}\in H$ such that
\begin{equation}\label{eq19}
z^{k}\rightharpoonup z^{\ast} \quad as \quad k\rightarrow +\infty,
\end{equation}
\begin{equation}\label{eq20}
w^{k} \rightharpoonup J_{\gamma B}z^{\ast} \quad as \quad k\rightarrow +\infty,
\end{equation}
\begin{equation}\label{eq21}
u^{k}\rightharpoonup J_{\gamma B}z^{\ast}=J_{\gamma A}(2J_{\gamma B}z^{\ast}-z^{\ast}-\gamma CJ_{\gamma B}z^{\ast})\quad as \quad k\rightarrow +\infty,
\end{equation}
\begin{equation}\label{eq22}
z^{k}-z^{k-1}\rightarrow 0 \quad as \quad k\rightarrow +\infty,
\end{equation}
\begin{equation}\label{eq23}
w^{k}-u^{k}\rightarrow 0 \quad as \quad k\rightarrow +\infty,
\end{equation}
\begin{equation}\label{eq24}
Cw^{k}\rightarrow CJ_{\gamma B}z^{\ast} \quad as \quad k\rightarrow +\infty.
\end{equation}

{\rm (i)}\, From (\ref{eq4}), we have $\gamma L_{3}x_{3}^{k}=z^{k}+\alpha_{k}(z^{k}-z^{k-1})-w^{k}+\gamma b$. By using (\ref{eq19}), (\ref{eq20}) and (\ref{eq22}), we get
\begin{equation}\label{eq25}
L_{3}x_{3}^{k}\rightharpoonup b+\frac{1}{\gamma}(z^{\ast}-J_{\gamma B}z^{\ast}).
\end{equation}

According to (\ref{eq5}), we have $L_{1}x_{1}^{k+1}=Cw^{k}$. Together with (\ref{eq24}), we obtain
\begin{equation}\label{eq26}
L_{1}x_{1}^{k+1} \rightarrow CJ_{\gamma B}z^{\ast}.
\end{equation}

From (\ref{eq8}), we know that $\gamma L_{2}x_{2}^{k+1}=2w^{k}-z^{k}-\alpha_{k}(z^{k}-z^{k-1})-\gamma L_{1}x_{1}^{k+1}-u^{k}$. By using (\ref{eq19}), (\ref{eq20}), (\ref{eq22}), (\ref{eq23}) and (\ref{eq26}), we get
\begin{equation}\label{eq27}
L_{2}x_{2}^{k+1} \rightharpoonup\frac{1}{\gamma}(J_{\gamma B}z^{\ast}-z^{\ast})-CJ_{\gamma B}z^{\ast}.
\end{equation}

Since for any $i\in \{1,2,3\}$, $\|L_{i}x_{i}\|\geq\theta_{i}\|x_{i}\|$, $\forall x_{i}\in H_{i}$. Then there exist $x_{1}^{\ast},x_{2}^{\ast}$, and $x_{3}^{\ast}$ such that
\begin{equation}\label{eq28}
L_{1}x_{1}^{\ast}=CJ_{\gamma B}z^{\ast}, \   L_{2}x_{2}^{\ast} = \frac{1}{\gamma}(J_{\gamma B}z^{\ast}-z^{\ast})-CJ_{\gamma B}z^{\ast}, \  \textrm{and} \   L_{3}x_{3}^{\ast}=b+\frac{1}{\gamma}(z^{\ast}-J_{\gamma B}z^{\ast}).
\end{equation}

According to (\ref{eq25}), (\ref{eq26}) and (\ref{eq27}), we get
\begin{equation}\label{eq29}
x_{1}^{k+1} \rightarrow x_{1}^{\ast}, \   x_{2}^{k+1} \rightharpoonup x_{2}^{\ast}, \  \textrm{and} \   x_{3}^{k} \rightharpoonup x_{3}^{\ast},
\end{equation}
and
\begin{equation}\label{eq30}
L_{1}x_{1}^{\ast}+L_{2}x_{2}^{\ast}+L_{3}x_{3}^{\ast}=b.
\end{equation}

Let $w^{\ast}=J_{\gamma B}z^{\ast}$, then
\begin{equation}\label{eq31}
z^{\ast}-w^{\ast}\in \gamma Bw^{\ast}.
\end{equation}
Thus $L_{3}x_{3}^{\ast}=b+\frac{1}{\gamma}(z^{\ast}-J_{\gamma B}z^{\ast})\in b+Bw^{\ast}$, which is equivalent to
\begin{equation}\label{eq32}
0\in \partial f_{3}(x_{3}^{\ast})-L_{3}^{\ast}w^{\ast}.
\end{equation}

According to $L_{1}x_{1}^{\ast}=CJ_{\gamma B}z^{\ast}=Cw^{\ast}=L_{1}\nabla f_{1}^{\ast}(L_{1}^{\ast}w^{\ast})$ which is
\begin{equation}\label{eq33}
0\in \partial f_{1}(x_{1}^{\ast})-L_{1}^{\ast}w^{\ast}.
\end{equation}

Again from (\ref{eq21}), $w^{\ast}=J_{\gamma A}(2w^{\ast}-z^{\ast}-\gamma Cw^{\ast})$, we have $\frac{1}{\gamma}(w^{\ast}-z^{\ast})-Cw^{\ast}\in Aw^{\ast}$. Since $L_{2}x_{2}^{\ast} = \frac{1}{\gamma}(J_{\gamma B}z^{\ast}-z^{\ast})-CJ_{\gamma B}z^{\ast}=\frac{1}{\gamma}(w^{\ast}-z^{\ast})-Cw^{\ast}$, we get
\begin{equation}\label{eq34}
0\in \partial f_{2}(x_{2}^{\ast})-L_{2}^{\ast}w^{\ast}.
\end{equation}

According to (\ref{eq30}), (\ref{eq32}), (\ref{eq33}) and (\ref{eq34}), we prove that point pair $(x_{1}^{\ast},x_{2}^{\ast},x_{3}^{\ast},w^{\ast})$ satisfies optimality condition (\ref{eq3.3}), that is, point pair $(x_{1}^{\ast},x_{2}^{\ast},x_{3}^{\ast},w^{\ast})$ is saddle point of Lagrangian function (\ref{eq3.1}).

{\rm (ii)} We can get it directly from $w^{k}\rightharpoonup J_{\gamma B}z^{\ast}$ and $w^{\ast}=J_{\gamma B}z^{\ast}$.

{\rm (iii)} From (\ref{eq9}) we have $L_{1}x_{1}^{k+1}+L_{2}x_{2}^{k+1}+L_{3}x_{3}^{k}= \frac{1}{\gamma}(w^{k}-u^{k})+b$, and then use (\ref{eq23}).

{\rm (iv)} Suppose that one of the following conditions holds:

\quad\quad (a) $f_{1}^{\ast}$ is uniformly convex on every nonempty bounded subset of $dom f_{1}^{\ast}$;

\quad\quad (b) $f_{2}^{\ast}$ is uniformly convex on every nonempty bounded subset of $dom f_{2}^{\ast}$;

\quad\quad (c) $f_{3}^{\ast}$ is uniformly convex on every nonempty bounded subset of $dom f_{3}^{\ast}$.

Assume that for any $i\in\{1,2,3\}$, $\|L_{i}^{\ast}x\|\geq\beta_{i}\|x\|,\forall x\in H$, for some $\beta_{i}>0$.

(a) Suppose that $f_{1}^{\ast}$ is uniformly convex. Let $x_{1}\in H,x_{2}\in H$, then there exists an nondecreasing function $\phi_{f_{1}^{\ast}}:[0,+\infty)\rightarrow[0,+\infty)$ that vanishes only at 0 such that
\begin{align}\label{eq35}
\langle x_{1}-x_{2},L_{1}\nabla f_{1}^{\ast}(L_{1}^{\ast}x_{1})-L_{1}\nabla f_{1}^{\ast}(L_{1}^{\ast}x_{2})\rangle & =\langle L_{1}^{\ast}x_{1}-L_{1}^{\ast}x_{2},\nabla f_{1}^{\ast}(L_{1}^{\ast}x_{1})-\nabla f_{1}^{\ast}(L_{1}^{\ast}x_{2})\rangle  \nonumber\\
& \geq \phi_{f_{1}^{\ast}}(\|L_{1}^{\ast}x_{1}-L_{1}^{\ast}x_{2}\|) \nonumber \\
& \geq \phi_{f_{1}^{\ast}}(\beta_{1}\|x_{1}-x_{2}\|),
\end{align}
which implies that $C=L_{1}\circ\nabla f_{1}^{\ast}\circ L_{1}^{\ast}$ is uniformly monotone.

Similarly, we can prove that $A=L_{2}\circ\partial f_{2}^{\ast}\circ L_{2}^{\ast}$ and $B=L_{3}\circ\partial f_{3}^{\ast}\circ L_{3}^{\ast}$ are also uniformly monotone under the conditions of (b) or (c).

{\rm (v)} We know that $f_{i}(i=1,2,3)$ are lower semi-continuous, therefore we have
\begin{align}\label{eq36}
\liminf_{k\rightarrow+\infty}(f_{1}(x_{1}^{k+1})+f_{2}(x_{2}^{k+1})+ f_{3}(x_{3}^{k})) & \geq \liminf_{k\rightarrow+\infty}f_{1}(x_{1}^{k+1}) +\liminf_{k\rightarrow+\infty}f_{2}(x_{2}^{k+1}) +\liminf_{k\rightarrow+\infty}f_{3}(x_{3}^{k}),\nonumber\\
& \geq f_{1}(x_{1}^{\ast})+f_{2}(x_{2}^{\ast})+ f_{3}(x_{3}^{\ast})=v(P).
\end{align}

On the other hand, since $L_{1}^{\ast}w^{k}\in\partial f_{1}(x_{1}^{k+1})$, $L_{2}^{\ast}u^{k}\in\partial f_{2}(x_{2}^{k+1})$ and $L_{3}^{\ast}w^{k}\in\partial f_{3}(x_{3}^{k})$, we get
\begin{equation}\label{eq37}
f_{1}(x_{1}^{\ast})\geq f_{1}(x_{1}^{k+1})+\langle x_{1}^{\ast}-x_{1}^{k+1},L_{1}^{\ast}w^{k}\rangle ,
\end{equation}
\begin{equation}\label{eq38}
f_{2}(x_{2}^{\ast})\geq f_{2}(x_{2}^{k+1})+\langle x_{2}^{\ast}-x_{2}^{k+1},L_{2}^{\ast}u^{k}\rangle,
\end{equation}
\begin{equation}\label{eq39}
f_{3}(x_{3}^{\ast})\geq f_{3}(x_{3}^{k})+\langle x_{3}^{\ast}-x_{3}^{k},L_{3}^{\ast}w^{k}\rangle.
\end{equation}
Adding (\ref{eq37})-(\ref{eq39}), we obtain
\begin{align}\label{eq40}
v(P)\geq f_{1}(x_{1}^{k+1})&+ f_{2}(x_{2}^{k+1})+ f_{3}(x_{3}^{k})\nonumber\\
& +\langle b-L_{1}x_{1}^{k+1}-L_{2}x_{2}^{k+1}-L_{3}x_{3}^{k}, w^{k}\rangle+\langle L_{2}x_{2}^{\ast}-L_{2}x_{2}^{k+1},u^{k}-w^{k}\rangle.
\end{align}

Again from (i), (ii), (iii) and (\ref{eq23}), we have
\begin{equation}\label{eq41}
\limsup_{k\rightarrow+\infty}(f_{1}(x_{1}^{k+1})+f_{2}(x_{2}^{k+1})+ f_{3}(x_{3}^{k}))\leq v(P).
\end{equation}
Combined with (\ref{eq36}) and (\ref{eq41}), we complete the first part of the Theorem \ref{theo3} (v).

From $L_{1}^{\ast}w^{k}\in\partial f_{1}(x_{1}^{k+1})$, $L_{2}^{\ast}u^{k}\in\partial f_{2}(x_{2}^{k+1})$ and $L_{3}^{\ast}w^{k}\in\partial f_{3}(x_{3}^{k})$, then
\begin{equation}\label{eq42}
f_{1}(x_{1}^{k+1})+f_{1}^{\ast}(L_{1}^{\ast}w^{k})=\langle x_{1}^{k+1},L_{1}^{\ast}w^{k}\rangle,
\end{equation}
\begin{equation}\label{eq43}
f_{2}(x_{2}^{k+1})+f_{2}^{\ast}(L_{2}^{\ast}u^{k})=\langle x_{2}^{k+1},L_{2}^{\ast}u^{k}\rangle,
\end{equation}
\begin{equation}\label{eq44}
f_{3}(x_{3}^{k})+f_{3}^{\ast}(L_{3}^{\ast}w^{k})=\langle x_{3}^{k},L_{3}^{\ast}w^{k}\rangle.
\end{equation}
Adding (\ref{eq42})-(\ref{eq44}), we obtain
\begin{align}\label{eq45}
f_{1}(x_{1}^{k+1})+f_{2}(x_{2}^{k+1})+ f_{3}(x_{3}^{k})= & -f_{1}^{\ast}(L_{1}^{\ast}w^{k})  -f_{2}^{\ast}(L_{2}^{\ast}u^{k})- f_{3}^{\ast}(L_{3}^{\ast}w^{k})\nonumber\\
& +\langle L_{1}x_{1}^{k+1}+L_{2}x_{2}^{k+1}+L_{3}x_{3}^{k}, w^{k}\rangle + \langle L_{2}x_{2}^{k+1},u^{k}-w^{k}\rangle.
\end{align}

Finally, taking into account (i), (iii), (\ref{eq23}) and the first part of Theorem \ref{theo3} (v), we get
\begin{align}\label{eq46}
\lim_{k\rightarrow+\infty}(-f_{1}^{\ast}(L_{1}^{\ast}w^{k}) -f_{2}^{\ast}(L_{2}^{\ast}u^{k})-f_{3}^{\ast}(L_{3}^{\ast}w^{k}) +\langle w^{k},b\rangle)=v(D)=v(P).
\end{align}
This completes the proof.
\end{proof}

\begin{theorem}\label{theo4}
Suppose that the assumptions (A1)-(A3) are valid. Let $\{(x_{1}^{k},x_{2}^{k},x_{3}^{k},w^{k})\}$ be the sequence generated by Algorithm \ref{algo2}. Let $\gamma\in(0,2\beta\bar{\varepsilon})$, where $\bar{\varepsilon}\in(0,1)$ and $\beta=\mu/\|L_{1}\|^{2}$. Let $0\leq\alpha_{k}\leq\alpha<1$ and $0<\underline{\lambda}\leq \lambda_{k}\bar{\alpha}\leq\overline{\lambda}<1$, where $\bar{\alpha}=\frac{1}{2-\bar{\varepsilon}}$. Let $\sum_{k=1}^{+\infty}\alpha_{k+1}\|p^{k}-\gamma\lambda_{k}(L_{1}x_{1}^{k+1} + L_{2}x_{2}^{k+1} + L_{3}x_{3}^{k} -b)\|^{2}<+\infty$. Then there exists a point pair $(x_{1}^{\ast},x_{2}^{\ast},x_{3}^{\ast},w^{\ast})$, which is the saddle point of the Lagrange function (\ref{eq3.1}) such that the following hold:

\rm(i) $\{(x_{1}^{k+1},x_{2}^{k+1},x_{3}^{k+1})\}_{k\geq1}$ converges weakly to $(x_{1}^{\ast},x_{2}^{\ast},x_{3}^{\ast})$. In particular, $\{x_{1}^{k+1}\}_{k\geq1}$ converges strongly to $x_{1}^{\ast}$;

(ii) $\{w^{k+1}\}_{k\geq1}$ converges weakly to $w^{\ast}$;

(iii) $\{L_{1}x_{1}^{k+1}+L_{2}x_{2}^{k+1}+L_{3}x_{3}^{k}\}_{k\geq2}$ converges strongly to $b$;

(iv) Suppose that one of the following conditions hold:

\quad\quad (a) $f_{1}^{\ast}$ is uniformly convex on every nonempty bounded subset of $dom f_{1}^{\ast}$;

\quad\quad (b) $f_{2}^{\ast}$ is uniformly convex on every nonempty bounded subset of $dom f_{2}^{\ast}$;

\quad\quad (c) $f_{3}^{\ast}$ is uniformly convex on every nonempty bounded subset of $dom f_{3}^{\ast}$;

then $\{w^{k+1}\}_{k\geq1}$ converges strongly to the unique optimal solution of $(D)$;

(v) $\lim_{k\rightarrow+\infty}(f_{1}(x_{1}^{k+1})+f_{2}(x_{2}^{k+1})+f_{3}(x_{3}^{k})) =v(P)=v(D)=\lim_{k\rightarrow+\infty}(-f_{1}^{\ast}(L_{1}^{\ast}w^{k}) -f_{2}^{\ast}(L_{2}^{\ast}u^{k})-f_{3}^{\ast}(L_{3}^{\ast}w^{k}) +\langle w^{k},b\rangle)$, where $u^{k}$ is defined as (\ref{eq1}).
\end{theorem}

\vskip 2mm

\begin{proof}
The proof of Theorem \ref{theo4} is similar to Theorem \ref{theo3}, so we omit it here.
\end{proof}
\vskip 2mm

\begin{remark}
Notice that, in finite-dimensional case, the assumption on $L_{i}$ $(i=1,2,3)$ in Theorem \ref{theo3} and Theorem \ref{theo4} means that $L_{i}$ are matrices with full column rank and all weak convergences in Theorem \ref{theo3} and Theorem \ref{theo4} are strong convergences.
\end{remark}
\vskip 2mm
\begin{remark}
In comparison with the other three-block ADMM, such as (\ref{3-block-ADM-G}) and (\ref{3-block-sPADMM}). We prove the weak and strong convergence of the iteration sequences generated by Algorithm \ref{algo2}. However, the strong convergence of three-block ADMM (\ref{3-block-ADM-G}) and (\ref{3-block-sPADMM}) are only proved in finite-dimensional Hilbert spaces. It's not clear whether they still have strong convergence in infinite-dimensional Hilbert spaces. It is well-known that the weak and strong convergence is not equivalent to each other in infinite-dimensional Hilbert spaces.
\end{remark}
\vskip 2mm

In the following, we present several particular cases of the proposed Algorithm \ref{algo2}.

Let $\alpha_{k}=0$ in Algorithm \ref{algo2}, then we get the relaxed three-block alternating minimization algorithm (R-AMA)
\begin{equation}
   \left\{
\begin{aligned}\label{3-block-R-AMA}
 & x_{1}^{k+1}   = \arg\min_{x_{1}} \{ f_{1}(x_{1})-\langle w^k,L_{1}x_{1}\rangle\},\\
 & x_{2}^{k+1}  = \arg\min_{x_{2}} \{ f_{2}(x_{2}) - \langle w^k , L_{2}x_{2} \rangle + \frac{\gamma}{2} \| L_{1}x_{1}^{k+1} + L_{2}x_{2} + L_{3}x_{3}^{k} -b \|^2  \},\\
 & x_{3}^{k+1}  = \arg\min_{x_{3}} \{ f_{3}(x_{3}) -\langle w^k, L_{3}x_{3} \rangle + \frac{\gamma}{2} \|L_{3}(x_{3}-x_{3}^{k})+\lambda_{k}(L_{1}x_{1}^{k+1} + L_{2}x_{2}^{k+1} + L_{3}x_{3}^{k} -b) \|^2\},\\
 & w^{k+1}  = w^k - \gamma(L_{3}(x_{3}^{k+1}-x_{3}^{k}) +\lambda_{k}(L_{1}x_{1}^{k+1} + L_{2}x_{2}^{k+1} + L_{3}x_{3}^{k} -b)).
\end{aligned}
\right.
\end{equation}
Further, let $\lambda_{k} = 1$ in (\ref{3-block-R-AMA}), we recover the three-block AMA proposed by Davis and Yin \cite{davis2015}.

In Algorithm \ref{algo2}, when $f_{1}$ and $x_{1}$ vanish, the iterative sequences of Algorithm \ref{algo2} becomes for every $k\geq 1$,
\begin{equation}
   \left\{
\begin{aligned}\label{2-block-iADMM}
 & x_{2}^{k+1}  = \arg\min_{x_{2}} \{ f_{2}(x_{2}) - \langle w^k , L_{2}x_{2} \rangle + \frac{\gamma}{2} \| L_{2}x_{2} + L_{3}x_{3}^k -b \|^2  \},\\
 & x_{3}^{k+1}  = \arg\min_{x_{3}} \{ f_{3}(x_{3}) -\langle w^k+\alpha_{k+1}p^{k}, L_{3}x_{3} \rangle + \frac{\gamma}{2} \|L_{3}(x_{3}-x_{3}^{k}) +(1+\alpha_{k+1})\lambda_{k}(L_{2}x_{2}^{k+1} + L_{3}x_{3}^{k} -b) \|^2\},\\
 & w^{k+1}  = w^k + \alpha_{k+1}p^{k} - \gamma(L_{3}(x_{3}^{k+1}-x_{3}^{k}) +(1+\alpha_{k+1})\lambda_{k}(L_{2}x_{2}^{k+1} + L_{3}x_{3}^{k} -b)),\\
 & p^{k+1} =\alpha_{k+1}(p^{k}-\gamma\lambda_{k}(L_{2}x_{2}^{k+1} + L_{3}x_{3}^{k} -b)),
\end{aligned}
\right.
\end{equation}
which is the two-block inertial ADMM proposed in \cite{Yang2020}. Moreover, when $f_{2},x_{2}$ and $f_{3},x_{3}$ vanish respectively, and $\alpha_{k}=0$, then the Algorithm \ref{algo2} reduces to the following two different relaxed alternating minimization algorithms
\begin{equation}
   \left\{
\begin{aligned}\label{2-block-R-AMA2}
 & x_{1}^{k+1}   = \arg\min_{x_{1}} \{ f_{1}(x_{1})-\langle w^k,L_{1}x_{1}\rangle\},\\
 & x_{3}^{k+1}  = \arg\min_{x_{3}} \{ f_{3}(x_{3}) - \langle w^k , L_{3}x_{3} \rangle + \frac{\gamma}{2} \| L_{3}(x_{3}-x^{k}_{3}) +\lambda_{k}(L_{1}x_{1}^{k+1} + L_{3}x_{3}^{k} -b) \|^2  \},\\
 & w^{k+1}  = w^k - \gamma(L_{3}(x_{3}^{k+1}-x^{k}_{3}) +\lambda_{k}(L_{1}x_{1}^{k+1} + L_{3}x_{3}^{k} -b)),
\end{aligned}
\right.
\end{equation}
and
\begin{equation}
   \left\{
\begin{aligned}\label{2-block-R-AMA}
 & x_{1}^{k+1}   = \arg\min_{x_{1}} \{ f_{1}(x_{1})-\langle w^k,L_{1}x_{1}\rangle\},\\
 & x_{2}^{k+1}  = \arg\min_{x_{2}} \{ f_{2}(x_{2}) - \langle w^k , L_{2}x_{2} \rangle + \frac{\gamma}{2} \| L_{1}x_{1}^{k+1} + L_{2}x_{2} -b \|^2  \},\\
 & w^{k+1}  = w^k - \gamma\lambda_{k}(L_{1}x_{1}^{k+1} + L_{2}x_{2}^{k+1} -b).
\end{aligned}
\right.
\end{equation}
Let $\lambda_{k} = 1$, then (\ref{2-block-R-AMA2}) and (\ref{2-block-R-AMA}) are reduced to the alternating minimization algorithm (AMA) proposed by Tseng \cite{Tseng1991SIAM}.

\vskip 2mm

\section{Numerical experiments}
\label{numer_test}
In this section, we carry out simulation experiments and compare the proposed algorithm (Algorithm \ref{algo2}) and its by-product relaxed alternative minimization algorithm (R-AMA (\ref{3-block-R-AMA})) with other state-of-the-art algorithms include the three-block ADMM (\ref{3-block-ADMM}) \cite{Cai2017COA}, the ADM-G (\ref{3-block-ADM-G}) \cite{He2012SJOO}, the sPADMM (\ref{3-block-sPADMM}) \cite{LI2015APJOOR} and three-block alternative minimization algorithm (AMA (\ref{3-block-AMA})) proposed by Davis and Yin \cite{davis2015} on the stable principal component pursuit (SPCP). All the experiments are conducted on a 64-bit Windows 10 operating system with an Intel(R) Core(TM) i5-7200U CPU and 8GB memory. All the codes are tested in MATLAB R2016a.

\subsection{Stable principal component pursuit (SPCP)}

The purpose of the stable principal component pursuit \cite{Zhou2010ISIT} is to recover the low-rank matrix from the high dimensional data matrix with sparse error and small noise. This problem is a special case of (\ref{problem1}), which can be formulated as:
\begin{equation}\label{SPCP}
\begin{aligned}
 \min_{L,S,Z}\, & \beta_{1}\|L\|_{\ast} + \beta_{2}\|S\|_{1} + \frac{1}{2}\|Z\|^{2}_{F}  \\
\textrm{s.t.}\, & L + S + Z = b,
\end{aligned}
\end{equation}
Where $\|L\|_{\ast}$ is defined as the sum of all singular values of the matrix $L$, $\|S\|_{1}$ is the $l_{1}$ norm of the matrix $S$ and $\|Z\|_{F}$ is the Frobenius norm of the matrix $Z$; $b$ is a given damaged data matrix, $L$, $S$ and $Z$ are a low rank, sparse and noise components of $b$, respectively. We conduct numerical experiments with the generated simulation data to show the effectiveness of the proposed algorithm. The generation of simulation data is similar to \cite{Candes2009ACM}. The observed damaged data matrix $b$ is generated as follows.
The low rank matrix $L^{\ast}$ is generated by $L^{\ast} = L_{1}L_{2}^{T}$, where $L_{1} = randn(m, r)$ and $L_{2}= randn(m, r)$ are two independently generated random matrices of $m\times r$ ($r<m$ is rank of matrix $L^{\ast}$) scale. The $S^{\ast}$ is a sparse matrix with non-zero elements uniformly distributed and values uniformly distributed between [-500,500]. The $Z^{\ast}$ is a matrix with Gaussian noise whose mean value is 0 and standard deviation is $10^{-5}$. Finally, we set $b=L^{\ast}+S^{\ast}+Z^{\ast}$.

We put the actual problem (\ref{SPCP}) into Algorithm \ref{algo2}. Let $x_{1}:=Z,x_{2}:=L$ and $x_{3}:=S$, it is obvious that problem (\ref{SPCP}) is a special case of model (\ref{problem1}). Accordingly, $f_{1}(x_{1}):=\frac{1}{2}\|Z\|_{F}^{2}$, $f_{2}(x_{2}):=\beta_{1}\|L\|_{\ast}$ and $f_{3}(x_{3}):=\beta_{2}\|S\|_{1}$, coefficient matrixes $L_{1}=L_{2}=L_{3}:=I$, where $I$ is the identity operator.

The following is the detailed calculation process of the problem (\ref{SPCP}) executing Algorithm \ref{algo2}.

\rm1. $Z$-subproblem in Algorithm \ref{algo2}:
\begin{align}
Z^{k+1}&=\arg\min_{Z}\{\frac{1}{2}\|Z\|^{2}-\langle w^{k},Z\rangle\} \nonumber\\
&=w^{k}.\nonumber
\end{align}

\rm2. $L$-subproblem in Algorithm \ref{algo2}:
\begin{align}
L^{k+1}&=\arg\min_{L}\{\beta_{1}\|L\|_{\ast}-\langle w^k , L \rangle + \frac{\gamma}{2} \|Z^{k+1} + L  + S^{k} -b \|^2  \} \nonumber\\
&=\arg\min_{L}\{\beta_{1}\|L\|_{\ast}+ \frac{\gamma}{2} \|Z^{k+1} + L  + S^{k} -b -\frac{1}{\gamma}w^{k}\|^2  \} \nonumber\\
&=prox_{\frac{\beta_{1}}{\gamma}\|\cdot\|_{\ast}}(b+\frac{1}{\gamma}w^{k}-Z^{k+1}-S^{k}),\nonumber
\end{align}
where $prox_{c\|\cdot\|_{\ast}}(\cdot)$ is the proximal function \cite{Caijf2010SIAM} of the function $c\|\cdot\|_{\ast}$ with a constant $c>0$. For any matrix $L\in R^{m\times n}$ with $rank(L)=r$, let its singular value decomposition be $L=Udiag(\{\sigma_{i}\}_{1\leq i\leq r})V^{T}$, where $U\in R^{m\times r}$, $V\in R^{n\times r}$, then $prox_{c\|\cdot\|_{\ast}}(L)=Udiag(max\{\{\sigma_{i}\}_{1\leq i\leq r}-c,0\})V^{T}$.

\rm3. $S$-subproblem in Algorithm \ref{algo2}:
\begin{align}
S^{k+1} & = \arg\min_{S} \{ \beta_{2}\|S\|_{1} -\langle w^k+\alpha_{k+1}p^{k}, S \rangle + \frac{\gamma}{2} \|(S-S^{k}) +(1+\alpha_{k+1})\lambda_{k}(Z^{k+1} + L^{k+1} + S^{k} -b) \|^2\} \nonumber\\
&=\arg\min_{S} \{ \beta_{2}\|S\|_{1}+ \frac{\gamma}{2}\|(S-S^{k}) +(1+\alpha_{k+1})\lambda_{k}(Z^{k+1} + L^{k+1} + S^{k} -b)-\frac{1}{\gamma}(w^k+\alpha_{k+1}p^{k}) \|^2\} \nonumber\\
&=prox_{\frac{\beta_{2}}{\gamma}\|\cdot\|_{1}}(S^{k}+\frac{1}{\gamma}(w^k+\alpha_{k+1}p^{k})- (1+\alpha_{k+1})\lambda_{k}(Z^{k+1} + L^{k+1} + S^{k} -b)).\nonumber
\end{align}
where $prox_{c\|\cdot\|_{1}}(S) = sign(S).*max(abs(S)-c,0)$.

\rm4. Update of Lagrange multiplier $w$ in Algorithm \ref{algo2}:
\begin{align}
w^{k+1}  = w^k + \alpha_{k+1}p^{k} - \gamma(S^{k+1}-S^{k} +(1+\alpha_{k+1})\lambda_{k}(Z^{k+1} + L^{k+1} + S^{k} -b)).\nonumber
\end{align}

\rm5. Update of variable $p$ in Algorithm \ref{algo2}:
\begin{align}
p^{k+1} =\alpha_{k+1}(p^{k}-\gamma\lambda_k(Z^{k+1}+L^{k+1}+S^{k}-b)).\nonumber
\end{align}

\vskip 2mm

\subsection{Parameters setting}

The specific setting of each parameter in the algorithm is given in this subsection. Let $\beta_{1}=0.05$ and $\beta_{2}=\beta_{1}/\sqrt{m}$. Let the relative error of $L$ and $S$ be the stopping criterion, i.e.,
\begin{align}\label{eq4.2.1}
rel\, L:=\frac{\|L^{k+1}-L^{k}\|_{F}}{\|L^{k}\|_{F}}&,\ rel\, S:=\frac{\|S^{k+1}-S^{k}\|_{F}}{\|S^{k}\|_{F}},\nonumber\\
\max(rel\,L &, rel\,S)\leq \varepsilon, \nonumber
\end{align}
where $\varepsilon$ is a small constant. We first conduct an numerical experiment to illustrate the relationship between the value of the penalty parameter $\gamma$ and the experimental results such as the number of iteration steps in the three-block AMA (\ref{3-block-AMA}) algorithm. In this experiment, we set $m=200$, $rank(L^{\ast})=0.05m$, $\|S^{\ast}\|_{0}=0.05m^{2}$, $\varepsilon=10^{-5}$, and the initial variables $(S^{1},w^{1})=(0,0)$.

\begin{table}[htbp]
\small
\centering
\caption{Numerical experimental results of three-block AMA (\ref{3-block-AMA}) algorithm under different penalty parameters $\gamma$($rel\ L^{\ast}$ and $rel\ S^{\ast}$ are defined as $\frac{\|L^{k}-L^{\ast}\|_{F}}{\|L^{\ast}\|_{F}}$ and $\frac{\|S^{k}-S^{\ast}\|_{F}}{\|S^{\ast}\|_{F}}$, respectively).}
\begin{tabular}{c|c|cccc}
\hline
\multirow{2}[1]{*}{Methods} & \multicolumn{5}{c}{$m = 200\ rank(L^{\ast})=0.05m\ \|S^{\ast}\|_{0}=0.05m^{2}\ \varepsilon=10^{-5}$}  \\
 \cline{2-6}
& $\gamma$ & $k$ & $rank(L^{k})$ & $rel\ L^{\ast}$ & $rel\ S^{\ast}$ \\
\hline
\hline
\multirow{9}[1]{*}{ AMA (\ref{3-block-AMA})} & 0.0005 & 37 & 10 & 3.2242e-4 & 2.3868e-5  \\
 \cline{2-6}
 & 0.005 & 224 & 10 & 2.8690e-4 & 1.4136e-5   \\
 \cline{2-6}
  & 0.05 & 2160 & 10 & 2.8738e-4 & 1.4145e-5   \\
 \cline{2-6}
  & 0.1 & 4311 & 10 & 2.8814e-4 & 1.4161e-5   \\
 \cline{2-6}
  & 0.5 & 21534 & 10 & 2.8719e-4 & 1.4140e-5   \\
  \cline{2-6}
  & 1 & 43066 & 10 & 2.8759e-4 & 1.4150e-5   \\
 \cline{2-6}
 & 1.2 & 51679 & 10 & 2.8765e-4 & 1.4151e-5   \\
 \cline{2-6}
 & 1.5 & 64598 & 10 & 2.8778e-4 & 1.4153e-5   \\
 \cline{2-6}
 & 1.8 & 77516 & 10 & 2.8853e-4 & 1.4168e-5   \\
 \hline
\end{tabular}\label{table1}
\end{table}
From Table \ref{table1}, we can see that when the value of the penalty parameter $\gamma$ is large, the iteration step $k$ of the three-block AMA (\ref{3-block-AMA}) algorithm is large. When $\gamma=0.0005$, the 3-block AMA (\ref{3-block-AMA}) algorithm has the fastest convergence speed. In the following experiments, we fix the $\gamma=0.0005$ and compare the effects of different relaxation parameters $\lambda_{k}$ on the numerical experimental results of the three-block R-AMA (\ref{3-block-R-AMA}) algorithm. The above experimental data is still used, and the $\lambda_{k}$ takes ten different values as $0.5$, $0.8$, $1$, $1.1$, $1.2$, $1.3$, $1.5$, $1.6$, $1.7$ and $1.8$, respectively.

\begin{table}[htbp]
\small
\centering
\caption{Numerical experimental results of three-block R-AMA (\ref{3-block-R-AMA}) algorithm under different relaxation parameters $\lambda_{k}$($rel\ L^{\ast}$ and $rel\ S^{\ast}$ are defined as $\frac{\|L^{k}-L^{\ast}\|_{F}}{\|L^{\ast}\|_{F}}$ and $\frac{\|S^{k}-S^{\ast}\|_{F}}{\|S^{\ast}\|_{F}}$, respectively).}
\begin{tabular}{c|c|c|cccc}
\hline
\multirow{2}[1]{*}{Methods} & \multicolumn{6}{c}{$m = 200\ rank(L^{\ast})=0.05m\ \|S^{\ast}\|_{0}=0.05m^{2}\ \varepsilon=10^{-5}$}  \\
 \cline{2-7}
& $\gamma$ & $\lambda_{k}$ & $k$ & $rank(L^{k})$ & $rel\ L^{\ast}$ & $rel\ S^{\ast}$ \\
\hline
\hline
\multirow{10}[1]{*}{ R-AMA (\ref{3-block-R-AMA})} & \multirow{10}[1]{*}{0.0005} & 0.5 & 69 & 10 & 3.2301e-4 &   2.3880e-5 \\
 \cline{3-7}
 &  & 0.8 & 46 & 10 & 3.2245e-4 & 2.3867e-5  \\
 \cline{3-7}
 &  & 1   & 37 & 10 & 3.2242e-4 & 2.3868e-5  \\
 \cline{3-7}
 &  & 1.1 & 34 & 10 & 3.2209e-4 & 2.3866e-5  \\
 \cline{3-7}
 &  & 1.2 & 31 & 10 & 3.2202e-4 & 2.3865e-5  \\
 \cline{3-7}
 &  & 1.3 & 29 & 10 & 3.2210e-4 & 2.3866e-5  \\
 \cline{3-7}
 &  & 1.5 & 27 & 10 & 3.2210e-4 & 2.3865e-5  \\
 \cline{3-7}
 &  & 1.6 & 33 & 10 & 2.8743e-4 & 1.4149e-5  \\
 \cline{3-7}
 &  & 1.7 & 36 & 10 & 2.8743e-4 & 1.4157e-5 \\
 \cline{3-7}
 &  & 1.8 & 43 & 10 & 2.8753e-4 & 1.4203e-5  \\
 \hline
\end{tabular}\label{table2}
\end{table}

From Table \ref{table2}, we can see that the relaxation parameter $\lambda_{k}$ can effectively improve the convergence speed of the AMA algorithm. When $\lambda_{k}\in[1.1,1.7]$, the relaxation parameter $\lambda_{k}$ can accelerate the three-block AMA (\ref{3-block-AMA}) algorithm, and the optimal acceleration effect is $\lambda_{k}=1.5$. In the following experiments we fix the relaxation parameter $\lambda_{k}=1.5$ of the relaxed three-block AMA (\ref{3-block-R-AMA}) algorithm. Subsequently, we compare the three-block ADMM (\ref{3-block-ADMM}), ADM-G (\ref{3-block-ADM-G}), sPADMM (\ref{3-block-sPADMM}), AMA (\ref{3-block-AMA}), R-AMA (\ref{3-block-R-AMA}) and Algorithm \ref{algo2} with different conditions. When $\gamma = 0.0005$, it satisfies the three-block ADMM (\ref{3-block-ADMM}), ADM-G (\ref{3-block-ADM-G}) and sPADMM (\ref{3-block-sPADMM}) restrictions on penalty parameters. Make the parameter $\theta = 0.99999$ in ADM-G (\ref{3-block-ADM-G}) and $\tau = 1.2$ in sPADMM (\ref{3-block-sPADMM}). We know that $\mu=1$, $L_{1}=I$, that is, $\beta=\mu/\|L_{1}\|^{2}=1$. And $\gamma\in(0,2\beta\bar{\varepsilon})$, so we make $\bar{\varepsilon}=0.00026$ and $\bar{\alpha}=\frac{1}{2-\bar{\varepsilon}}\approx0.5001$. We define their parameters in Table \ref{table3}.

\begin{table}[htbp]
\small
\centering
\caption{Parameters selection of the compared iterative algorithms.}
\begin{tabular}{c|c|c|c}
\hline
 Methods & $\gamma$  & $\lambda_{k}$ &  Inertial parameter $\alpha_{k}$ \\
 \hline
 \hline
 three-block ADMM (\ref{3-block-ADMM}) & \multirow{7}[1]{*}{$0.0005$} & None & None \\
 ADM-G (\ref{3-block-ADM-G}) &  & None & None \\
 sPADMM (\ref{3-block-sPADMM}) &  & None & None \\
 AMA (\ref{3-block-AMA})  & & 1 & None \\
 R-AMA (\ref{3-block-R-AMA})  & & 1.5 & None  \\
 Algorithm \ref{algo2}-1 &  & 1.25 & 0.15 \\
 Algorithm \ref{algo2}-2 & & 1.5 & $\min\{\frac{1}{k^{2}\|p^{k}-\gamma\lambda_{k}(L_{1}x_{1}^{k+1} + L_{2}x_{2}^{k+1} + L_{3}x_{3}^{k} -b)\|^{2}},0.005\}$ \\
\hline
\end{tabular}\label{table3}
\end{table}

\subsection{Results and discussions}

In order to make the experimental results more convincing, we conduct a number of numerical experiments. Let the order $m$ of the matrix be $200$, $400$ and $500$, respectively. The rank of low rank matrix $L^{\ast}$ and the sparsity of sparse matrix $S^{\ast}$ are also divided into two combinations: $rank(L^{\ast}) = 0.05m$ and $\|S^{\ast}\|_{0}=0.05m^{2}$, $rank(L^{\ast}) = 0.1m$ and $\|S^{\ast}\|_{0}=0.1m^{2}$.

We test the performance of the studied iterative algorithms including three-block ADMM (\ref{3-block-ADMM}), ADM-G (\ref{3-block-ADM-G}), sPADMM (\ref{3-block-sPADMM}), three-block AMA (\ref{3-block-AMA}), three-block R-AMA (\ref{3-block-R-AMA}), Algorithm \ref{algo2}-1 and Algorithm \ref{algo2}-2 with parameters selection in Table \ref{table3}. The results of numerical experiments are reported in Table \ref{table4}. Several indicators are listed here including the number of iteration steps, error accuracy, and running CPU time. From Table \ref{table4}, we can find that both the three-block R-AMA (\ref{3-block-R-AMA}) algorithm and the two relaxed inertial three-block AMA (Algorithm \ref{algo2}) algorithms with different conditions can accelerate the convergence speed of the three-block AMA (\ref{3-block-AMA}) algorithm, and their accuracy is higher. Table \ref{table4} also conveys a message: Inertia technology does not seem to be able to effectively accelerate the three-block AMA (\ref{3-block-AMA}) algorithm. The numerical performance of the two relaxed inertial three-block AMA (Algorithm \ref{algo2}) algorithms is almost the same or slightly worse than the three-block R-AMA (\ref{3-block-R-AMA}) algorithm. However, their performance is not as good as the three-block ADMM (\ref{3-block-ADMM}), ADM-G (\ref{3-block-ADM-G}) and sPADMM (\ref{3-block-sPADMM}). The iteration speed of the three-block ADMM (\ref{3-block-ADMM}) and sPADMM (\ref{3-block-sPADMM}) are almost the same, and they are faster than ADM-G (\ref{3-block-ADM-G}), which further proves that the direct promotion of the three-block ADMM (\ref{3-block-ADMM}) numerical experiment is better than other variants of ADMM.

\begin{table}[htbp]
\footnotesize
\newcommand{\tabincell}[2]{\begin{tabular}{@{}#1@{}}#2\end{tabular}}  
\centering
\caption{Comparison of numerical experimental results of three-block ADMM, ADM-G, sPADMM, AMA, R-AMA and Algorithm \ref{algo2}
($rel\ L^{\ast}$ and $rel\ S^{\ast}$ are defined as $\frac{\|L^{k}-L^{\ast}\|_{F}}{\|L^{\ast}\|_{F}}$ and $\frac{\|S^{k}-S^{\ast}\|_{F}}{\|S^{\ast}\|_{F}}$, respectively).}
\begin{tabular}{c|c|c|ccccc}
\hline
 & m & Methods  & $k$ & $rank(L^{k})$ & $rel\ L^{\ast}$ & $rel\ S^{\ast}$ & CPU \\
\hline
\hline
 \multirow{21}[1]{*}{\tabincell{c}{$rank(L^{\ast})=0.05m$ \\
 $\|S^{\ast}\|_{0}=0.05m^{2}$ \\ $\varepsilon=10^{-5}$ }} & \multirow{7}[1]{*}{200} & three-block ADMM (\ref{3-block-ADMM}) & 20 & 10 & 3.2221e-4 & 2.3867e-5 & 0.1946\\
\cline{3-8}
 & & ADM-G (\ref{3-block-ADM-G}) & 21 & 10 & 3.2216e-4 & 2.3858e-5 & 0.1934\\
\cline{3-8}
 & & sPADMM (\ref{3-block-sPADMM}) & 17 & 10 & 3.2219e-4 & 2.3868e-5 & 0.1740\\
\cline{3-8}
 & & AMA (\ref{3-block-AMA})  & 37 & 10 & 3.2242e-4 & 2.3868e-5 & 0.3450\\
\cline{3-8}
 & & R-AMA (\ref{3-block-R-AMA}) & 27 & 10 & 3.2210e-4 & 2.3865e-5 & 0.2379 \\
 \cline{3-8}
  & & Algorithm \ref{algo2}-1 & 28 & 10 & 3.2216e-4 & 2.3865e-5  & 0.2709\\
\cline{3-8}
  & & Algorithm \ref{algo2}-2 & 27 & 10 & 3.2214e-4 & 2.3866e-5  & 0.2507 \\
\cline{2-8}
 & \multirow{7}[1]{*}{400} &three-block ADMM (\ref{3-block-ADMM}) & 15 & 20 & 1.9399e-4 & 3.0007e-5 & 0.7093\\
\cline{3-8}
 & & ADM-G (\ref{3-block-ADM-G}) & 23 & 20 & 1.5296e-4 & 1.5123e-5 & 1.1701\\
\cline{3-8}
 & & sPADMM (\ref{3-block-sPADMM}) & 17 & 20 & 1.6064e-4 & 1.8351e-5 & 0.8746\\
\cline{3-8}
 & & AMA (\ref{3-block-AMA})  & 38 & 20 & 1.8482e-4 & 2.6547e-5 & 2.3205 \\
\cline{3-8}
 & & R-AMA (\ref{3-block-R-AMA}) & 33 & 20 & 1.6071e-4 & 1.8348e-5 &  1.6953\\
\cline{3-8}
  & & Algorithm \ref{algo2}-1 & 34 & 20 & 1.6064e-4 & 1.8348e-5  & 1.7018\\
\cline{3-8}
  & & Algorithm \ref{algo2}-2 & 33 & 20 & 1.6070e-4 & 1.8349e-5  & 1.7227 \\
\cline{2-8}
 & \multirow{7}[1]{*}{500} &three-block ADMM (\ref{3-block-ADMM}) & 17 & 25 & 1.3543e-4 & 2.0612e-5 & 1.5268\\
\cline{3-8}
 & & ADM-G (\ref{3-block-ADM-G}) & 23 & 25 & 1.1591e-4 & 1.0546e-5 & 2.0847\\
\cline{3-8}
 & & sPADMM (\ref{3-block-sPADMM}) & 18 & 25 & 1.1576e-4 & 1.0526e-5 & 1.5371\\
\cline{3-8}
  & & AMA (\ref{3-block-AMA})  & 42 & 25 & 1.2165e-4 & 1.3823e-5 & 3.5990\\
\cline{3-8}
  & & R-AMA (\ref{3-block-R-AMA}) & 35 & 25 & 1.1694e-4 & 1.1138e-5 & 2.9207 \\
\cline{3-8}
 & & Algorithm \ref{algo2}-1 & 35 & 25 & 1.1697e-4 & 1.1137e-5  & 3.0514\\
\cline{3-8}
 &  & Algorithm \ref{algo2}-2 & 35 & 25 & 1.1695e-4 & 1.1138e-5  & 2.9553 \\
 \hline
 \hline
 \multirow{21}[1]{*}{\tabincell{c}{$rank(L^{\ast})=0.1m$ \\
 $\|S^{\ast}\|_{0}=0.1m^{2}$ \\ $\varepsilon=10^{-5}$ }} & \multirow{7}[1]{*}{200} & three-block ADMM (\ref{3-block-ADMM}) & 19 & 20 & 4.1054e-4 & 3.0653e-5 & 0.1733\\
\cline{3-8}
 & & ADM-G (\ref{3-block-ADM-G}) & 27 & 20 & 4.0819e-4 & 3.0568e-5 & 0.2612\\
\cline{3-8}
 & & sPADMM (\ref{3-block-sPADMM}) & 18 & 20 & 4.0954e-4 & 3.0639e-5 & 0.2018\\
\cline{3-8}
 & & AMA (\ref{3-block-AMA})  & 48 & 20 & 3.3991e-4 & 1.8898e-5 & 0.4429\\
\cline{3-8}
 & & R-AMA (\ref{3-block-R-AMA}) & 36 & 20 & 3.3905e-4 & 1.8876e-5 &  0.3082\\
\cline{3-8}
  & & Algorithm \ref{algo2}-1 & 37 & 20 & 3.3883e-4 & 1.8872e-5  & 0.3320\\
\cline{3-8}
  & & Algorithm \ref{algo2}-2 & 36 & 20 & 3.3898e-4 & 1.88774e-5  & 0.3014 \\
\cline{2-8}
 & \multirow{7}[1]{*}{400} &three-block ADMM (\ref{3-block-ADMM}) & 28 & 40 & 2.1397e-4 & 2.2150e-5 & 1.4032\\
\cline{3-8}
 & & ADM-G (\ref{3-block-ADM-G}) & 32 & 40 & 2.1420e-4 & 2.1907e-5 & 1.8402\\
\cline{3-8}
 & & sPADMM (\ref{3-block-sPADMM}) & 34 & 40 & 1.7984e-4 & 1.5447e-5 & 1.6659\\
\cline{3-8}
 & & AMA (\ref{3-block-AMA})  & 60 & 40 & 2.2231e-4 & 2.3951e-5 & 3.3423\\
\cline{3-8}
 & & R-AMA (\ref{3-block-R-AMA}) & 52 & 40 & 1.7343e-4 & 1.3852e-5 & 2.6572 \\
\cline{3-8}
 & & Algorithm \ref{algo2}-1 & 53 & 40 & 1.7347e-4 & 1.3853e-5  & 2.8001\\
\cline{3-8}
  & & Algorithm \ref{algo2}-2 & 52 & 40 & 1.7336e-4 & 1.3850e-5  & 2.8206 \\
\cline{2-8}
 & \multirow{7}[1]{*}{500} &three-block ADMM (\ref{3-block-ADMM}) & 32 & 50 & 1.3999e-4 & 1.2946e-5 & 2.6338\\
\cline{3-8}
 & & ADM-G (\ref{3-block-ADM-G}) & 41 & 50 & 1.3394e-4 & 1.1220e-5 & 3.4820\\
\cline{3-8}
 & & sPADMM (\ref{3-block-sPADMM}) & 28 & 50 & 1.4008e-4 & 1.2967e-5 & 2.5707\\
\cline{3-8}
 & & AMA (\ref{3-block-AMA})  & 74 & 50 & 1.4346e-4 & 1.3376e-5 & 7.9132\\
\cline{3-8}
 & & R-AMA (\ref{3-block-R-AMA}) & 52 & 50 & 1.4302e-4 & 1.0980e-5 & 5.3290 \\
\cline{3-8}
 & & Algorithm \ref{algo2}-1 & 53 & 50 & 1.4303e-4 & 1.0981e-5  & 4.5382 \\
\cline{3-8}
 & & Algorithm \ref{algo2}-2 & 52 & 50 & 1.4201e-4 & 1.1018e-5 & 5.9416  \\
 \hline
\end{tabular}\label{table4}
\end{table}

\section{Conclusions}

The alternating direction method of multipliers (ADMM) and the alternating minimization algorithm (AMA) are two common splitting methods for solving separable convex programming with linear equality constraints. Recently, Davis and Yin \cite{davis2015} generalized the AMA to the case of three-block AMA (\ref{3-block-AMA}). In this paper, we proposed a relaxed inertial three-block AMA (Algorithm \ref{algo2}), which is derived from the inertial three-operator splitting algorithm \cite{Cui2019}. The obtained algorithm generalized and recovered some existing algorithms. In particular, we obtain a relaxed three-block AMA (\ref{3-block-R-AMA}). We analyze the convergence of the proposed algorithm in infinite-dimensional Hilbert spaces. Compared with other three-block ADMM, our convergence conclusions have not only weak convergence but also strong convergence. To demonstrate the efficiency and effectiveness of the proposed algorithm, we conduct numerical experiments on the stable principal component pursuit \cite{Zhou2010ISIT}. Numerical results showed that the relaxed three-block AMA (\ref{3-block-R-AMA}) performs better than the three-block AMA (\ref{3-block-AMA}) when the relaxation parameter belongs to $[1.1,1.7]$. We also observed that the performance of the relaxed inertial three-block AMA is similar to the relaxed three-block AMA. Our numerical results also confirmed that the limitations of the inertial accelerated ADMM pointed by Poon and Liang \cite{Liang2019NIPs}.

Recently, Bitterlich et al. \cite{Bitterlich2019JOTA} proposed a proximal AMA, which added proximal terms to the subproblem of the original AMA. Therefore, we would like to present the first open question:

\noindent \textbf{Question 1.}\, \textit{ Can we study the convergence of the following proximal three-block AMA (\ref{proximal-3-block-AMA})?}

\begin{equation}\label{proximal-3-block-AMA}
   \left\{
\begin{aligned}
x_{1}^{k+1} &  = \arg\min_{x_{1}} \left\{ f_{1}(x_{1}) - \langle w^k , L_{1}x_{1} \rangle   + \frac{1}{2}\| x_1 - x_{1}^{k} \|_{M_{1}^{k}} \right\}, \\
x_{2}^{k+1} & = \arg\min_{x_{2}} \left\{ f_{2}(x_{2}) - \langle w^k , L_{2}x_{2} \rangle + \frac{\gamma}{2} \| L_{1}x_{1}^{k+1} + L_{2}x_{2} + L_{3}x_{3}^k -b \|^2  + \frac{1}{2}\| x_2 - x_{2}^{k} \|_{M_{2}^{k}} \right\}, \\
x_{3}^{k+1} & = \arg\min_{x_{3}} \left\{ f_{3}(x_{3}) - \langle w^k , L_{3}x_{3} \rangle + \frac{\gamma}{2} \| L_{1}x_{1}^{k+1} + L_{2}x_{2}^{k+1} + L_{3}x_{3} -b \|^2 + \frac{1}{2}\| x_3 - x_{3}^{k} \|_{M_{3}^{k}} \right\} , \\
w^{k+1} & = w^k - \gamma(L_{1}x_{1}^{k+1} + L_{2}x_{2}^{k+1} + L_{3}x_{3}^{k+1} -b),
\end{aligned}
\right.
\end{equation}
where $\{M_{1}^{k}\}$, $\{M_{2}^{k}\}$, and $\{M_{3}^{k}\}$ are self-adjoint positive semidefinite operators.

As we know, the AMA is equivalent to the forward-backward splitting algorithm applied to the corresponding dual problem. In 2013, Raguet et al. \cite{Raguet-SIAM-2013} proposed a generalized forward-backward splitting algorithm for finding a zero of the sum of a cocoercive operator $B$ (See Definition \ref{def2}) and a finite sum of maximally monotone operators $\{A\}_{i=1}^{m}$, that is, find $x\in H$, such that $0\in Bx + \sum_{i=1}^{m}A_i x$. It is natural to employ the generalized forward-backward splitting algorithm to solve the dual of the following multi-block convex separable optimization problem.
\begin{equation}\label{multi-block}
\begin{aligned}
\min_{x_1, \cdots, x_m}\, & \sum_{i=1}^{m}f_{i}(x_i) \\
s.t. & L_1 x_1 + L_2 x_2 + \cdots + L_m x_m =b,
\end{aligned}
\end{equation}
where $\{f_i\}_{i=1}^{m}: H_i \rightarrow (-\infty, +\infty]$ are proper, lower semicontinuous convex functions, $\{L_i\}_{i=1}^{m}:H_i \rightarrow H$ are nonzero bounded linear operators, and $f_1$ is a strongly convex function.  Then, we raise the second open question:

\noindent \textbf{Question 2.}\, \textit{Can we obtain a primal-dual iteration scheme for solving (\ref{multi-block}) from the generalized forward-backward splitting algorithm?}

\section*{Funding}

This work was funded by the National Natural Science Foundations of China (12061045, 11661056, 11771347, 12031003).


\end{document}